\documentclass[]{article}   
\usepackage{amsmath,amsfonts,amsthm,amssymb,dsfont,mathtools}

\usepackage{amssymb}
\usepackage{graphicx}
\usepackage{enumerate}

\usepackage[normalem]{ulem}
\usepackage[dvipsnames,svgnames]{xcolor}
\usepackage{fancyvrb}
\usepackage{subfigure}
\usepackage{epstopdf} 
\theoremstyle{plain}
\newtheorem{theorem}{Theorem}[section]
\newtheorem{lemma}[theorem]{Lemma}
\newtheorem{proposition}[theorem]{Proposition}
\newtheorem{corollary}[theorem]{Corollary}

\theoremstyle{definition}
\newtheorem{definition}[theorem]{Definition}

\theoremstyle{remark}
\newtheorem{remark}[theorem]{Remark}

\def\CC{\hbox{C\kern -.58em {\raise .54ex \hbox{$\scriptscriptstyle |$}}
  \kern-.55em {\raise .53ex \hbox{$\scriptscriptstyle |$}} }}

\def\eg{{\sl \thinspace e.g.},\ }
\def\ie{{\sl \thinspace i.e.},\ }

\def\ZZ{{{\rm Z}\kern-.28em{\rm Z}}}
\def\RR{\mathop{{\rm I}\kern-.2em{\rm R}}\nolimits}
\def\un{\hbox{\rm1\kern-.28em\hbox{I}}}

\newcommand{\E}{\mathbb{E}}
\def\co#1#2{{\textstyle({#1\atop#2})}}

\newcommand{\T}{\mathbb{T}}

\newcommand{\spl}{\mathbb{S}}

\renewcommand{\leq}{\leqslant}
\renewcommand{\geq}{\geqslant}
\renewcommand{\phi}{\varphi}
\renewcommand{\Sigma}{\varSigma}
\newcommand{\id}{{\hbox{1\kern-.28em\hbox{I}}}}

\DeclareMathOperator{\on}{([a,b];\T)}\newcommand{\e}{\varepsilon}

\def\un{\hbox{\rm1\kern-.28em\hbox{I}}}

\title{Piecewise Extended Chebyshev  Spaces:\\ a numerical test for design
}
\author{
Carolina Vittoria Beccari$^1$, Giulio Casciola$^1$, Marie-Laurence Mazure$^2$}
\date{\it \small
$^1$ Department of Mathematics, University of Bologna, P.zza di Porta San Donato 5, \\
40126 Bologna, Italy\\
$^2$ Universit\'e Grenoble Alpes,  Laboratoire Jean Kuntzmann, 
CNRS, UMR 5224, \\
BP 53, F-38041 Grenoble 9, France \\
\rm carolina.beccari2@unibo.it, giulio.casciola@unibo.it,
mazure@imag.fr
}

\begin{document}
\thispagestyle{empty}
\maketitle
\thispagestyle{empty}
\pagestyle{plain}

\begin{abstract}

Given a number  of  Extended Chebyshev (EC) spaces on adjacent intervals, all of the same dimension, we join them via convenient connection matrices without increasing the dimension. The global space is called a Piecewise Extended Chebyshev (PEC) Space. In such a space one can count the total number of zeroes of any non-zero element, exactly as in each EC-section-space. When this number is bounded above in the global space the same way as in its section-spaces,  we say that it is an Extended Chebyshev Piecewise (ECP) space.  A thorough study of ECP-spaces has been developed in the last two decades in relation to  blossoms, with a view to design. In particular, extending a classical procedure for EC-spaces,  ECP-spaces were recently proved to all be obtained by means of piecewise generalised derivatives. This yields an interesting  constructive characterisation of ECP-spaces. Unfortunately, except for low dimensions and for very few adjacent intervals, this characterisation proved to be rather difficult to handle  in practice. To try to overcome this difficulty, in the present  article we show how to reinterpret the constructive characterisation as  a theoretical procedure to determine whether or not a given PEC-space is an ECP-space. This procedure is then translated into a numerical test, whose usefulness is illustrated by relevant examples.

\end{abstract}

\noindent
{\bf Keywords:} {Extended Chebyshev (piecewise) spaces,  connection matrices, Bernstein-type bases, 
(piecewise) generalised derivatives,  blossoms, geometric design}\par
\medskip
\noindent
{\bf AMS subject classification:} 65D05, 65D17

\section{Introduction}

 By their ability to ensure unisolvence of Hermite interpolation problems or, equivalently, by the bound on the number of zeroes of their non-zero elements, Extended Chebyshev spaces are known as the most natural generalisations of polynomial spaces, and for this reason they are old tools in Approximation Theory \cite{karlin, schu}. In that direction they are generally defined by means of generalised derivatives associated with systems of weight functions, which permits to extend to them various well-known notions of the polynomial framework, \eg generalised divided differences \cite{GM} and associated Newton-type decompositions, Taylor formul{\ae}, \dots \cite{GM, schu, lyche}.

 Initiated by H. Pottmann \cite{HP}, the theory of Chebyshevian blossoming has permitted a deeper understanding of Extended Chebyshev spaces and Chebyshevian splines (\ie splines with pieces taken from the same Extended Chebyshev space and with ordinary  continuity at the knots \cite{bister, tina}), while enhancing  their resemblance
with polynomial (spline) spaces in connection with geometric design. The present paper is not at all a paper on blossoms, but it would not exist without the fundamental contribution of these powerful and elegant tools. In any situation where blossoms arise, the major difficulty consists in proving their pseudoaffinity in each variable which  extends the well-known affinity in each variable of polynomial blossoms. Once this proven, the classical design algorithms are somehow inherent in  Chebyshevian  blossoms which also guarantee shape preservation of the resulting Bernstein bases \cite{optimal}.
As a  recent important progress arising from blossoms, let us mention the complete description of all possible systems of weight functions which can be associated with a given Extended Chebyshev  space on a closed bounded interval \cite{JAT11}.
 Moreover the geometrical nature of Chebyshevian blossoms makes them ideal tools to express geometric contact between parametric curves. This naturally produces blossoms for Chebyshevian splines with similar properties and consequences, \eg geometric design algorithms,  B-spline bases,  shape preservation.

Unlike polynomials, Extended Chebyshev spaces explicitly or implicitly involve  shape parameters and this explains  why they offer more possibilities  in the control of shapes of curves and surfaces. To take full advantage of the Chebyshevian framework,
it is useful to consider {\it Piecewise Chebyshevian splines}, that is, splines with pieces taken from different Extended Chebyshev spaces all of the same dimension, the continuity between consecutive pieces being controlled by connection matrices. These splines were first considered by P.J. Barry in \cite{barry}, see also \cite{CA99, JAT99, GM1} and  \cite{goodman, DM} for geometrically continuous polynomial splines. To be of interest for applications, and in particular for geometric design, such a spline space $\spl$ is expected to possess a B-spline basis -- in the usual sense of a normalised basis composed of minimally supported splines-- and this feature should be maintained after knot insertion.  As a matter of fact, this requirement was proved to be equivalent to the existence of blossoms in the space $\spl$ \cite{beyond, equiv, readyspline}. Moreover, for an efficient control of the shapes, the B-spline bases are additionally expected to  be totally positive. This property can  automatically be derived from the properties of blossoms, and in particular from their pseudoaffinity. This explains why the terminology {\it ``\,$\spl$ is good for design"} was adopted whenever the spline space $\spl$ possesses blossoms. It should be mentioned that the interest of Piecewise Chebyshevian spline spaces good for design is not limited to design: they also naturally produce  multiresolution analyses with associated piecewise Chebyshevian wavelets \cite{TML}, they permit  approximation by Schoenberg-type operators \cite{JAT13}, they have useful applications in Isogeometric Analysis \cite{carla}, \dots

In the present article we focus on the special case where all interior knots have zero multiplicities, the corresponding spline spaces being referred to as {\it Piecewise Extended Chebyshev spaces} (PEC). The  first motivation  to consider this case is that zero multiplicities can efficiently be used to strengthen the shape effects \cite{PJMLV}.    The second motivation lies in the fact that determining the class of all piecewise Chebyshevian spline spaces which are good for design amounts to determining the class of all PEC-spaces which are good for design. Indeed, it was recently proved that a piecewise Chebyshevian  spline space is good for design if and only if it is based on a PEC-space good for design, that is, possessing blossoms \cite{equiv, readyspline}.  It is known that a given PEC-space $\E$ which contains constants possesses blossoms if and only if the PEC-space obtained from $\E$ by differentiation is an {Extended Chebyshev piecewise space} (ECP) in the sense that the global bound on the number of zeroes is exactly the same as in each of its section-spaces. The presence of blossoms in $\E$ can also be characterised  by the existence of systems of weight functions associated with the section-spaces, relative to which the continuity conditions are expressed by identity matrices  \cite{JAT07}. Given that we know how to obtain all possible systems of weight functions associated with the section-spaces, this characterisation naturally  provides us with a procedure to determine whether or not a given PEC-space is good for design. Nevertheless, in general this procedure proves to be all the more difficult to carry out  in practice as zero multiplicities allow no freedom between consecutive sections.  This motivated the search for a an effective numerical procedure as a replacement, to which the present work is devoted.

The  paper is organised as follows. The necessary background  is presented in Section 2, with special insistence on Bernstein and Bernstein-like bases and their behaviour under possible piecewise generalised derivatives associated with piecewise weight functions according to a process similar to the non-piecewise case. In Section 3, these results are first reinterpreted as a theoretical test to answer the question: {\it is a given PEC-space an ECP-space?}, which is in turn  transformed into a numerical test.  What we actually test is: can we repeatedly diminish the dimension via piecewise generalised derivatives?   We illustrate this test by relevant examples  in Section 4, in particular with a view to design with shape parameters. We conclude the paper with some comments on both the usefulness and the limits of the numerical procedure.
\section{Background}

In this section we briefly survey  the main results on $(n+1)$-dimensional piecewise spaces obtained from $(n+1)$-dimensional section-spaces on adjacent intervals joined by connecting left/right derivatives at the interior knots by appropriate matrices. These results were proved in many earlier articles by the third author to which we refer the reader to, \eg \cite{NUMA05, ready, JAT07, howto} and other references therein. This survey is deliberately presented  in  a way to facilitate the next section.
\subsection{Piecewise spaces via connection matrices}

Throughout this article we consider  a fixed  interval $[a,b]$, $a<b$, and a fixed sequence $\T=(t_1, \ldots, t_q)$ of $q\geq 1$ knots interior to $[a,b]$, with
$$
t_0:=a<t_1<\dots<t_q<t_{q+1}:=b.
$$
We will deal with {\it piecewise functions on $\on$}, defined separately on each $[t_k^+, t_{k+1}^-]$. Given two such piecewise functions $F, G$ on $\on$, the equality $F=G$ means that $F(x)=G(x)$ for all $x\in[a,b]\setminus\{t_1, \ldots, t_q\}$, and that both $F(t_k^-)=G(t_k^-)$ and  $F(t_k^+)=G(t_k^+)$, for $k=1, \ldots, q$. This will be summarised by saying that $F(x^\e)=G(x^\e)$ for all $x\in[a,b]$,  $\e $ having the meaning of both $-, +$ if $x\in]a,b[$. One can similarly consider positive piecewise functions on $\on$, and so forth. Whenever necessary we will allow the integer  $q$ to be zero to come back to the non-piecewise situation.

Throughout the article, $D$ will stand for the (possibly left/right) ordinary differentiation, and $\un$ for the constant function $\un(x)=1$ for all $x$, on any interval.

With a view to defining piecewise spaces, take:
 \begin{enumerate}[--]
\item a sequence $\E_k$, $0\leq k\leq  q$, of  section-spaces: for each $k$, $\E_k\subset C^n([t_k,t_{k+1}])$ is an $(n+1)$-dimensional W-space on $[t_k,t_{k+1}]$ (\ie the Wronskian of a basis of $\E_k$ never vanishes on $[t_k,t_{k+1}]$, or any Taylor interpolation in $(n+1)$ data at any $x\in [t_k,t_{k+1}]$ has a unique solution  in $\E_k$); 
\item a  sequence $M_1, \ldots, M_q$ of connection matrices of order $(n+1)$ : each $M_k$ is  lower triangular with positive diagonal entries.
\end{enumerate}
These ingredients provide us with an $(n+1)$-dimensional {\it Piecewise W- space} (for short, PW-space) {\it on $\on$}, defined as the set $\E$ of all {\it piecewise functions} $F$ on $\on$  such that
\begin{enumerate}[\rm \arabic{enumi})]
\item for $k=0,\dots,q$,  there exists a function $F_k\in \E_k$ such that  $F$ coincides with $F_k$ on $[t_k^+,t_{k+1}^-]$;
\item for $k=1,\dots,q$, the following connection condition is fulfilled:
\begin{equation}
\label{connection}
\left(F(t_k^+), F'(t_k^+),\dots,F^{(n)}(t_k^+) \right)^T = M_k \left( F(t_k^-), F'(t_k^-),\dots,F^{(n)}(t_k^-) \right)^T.
\end{equation}
\end{enumerate}
This is the most natural and the largest framework to define piecewise spaces by connecting left/right derivatives at the interior knots. In such a space, the Wronskian of any basis never vanishes on $\on$.

At this stage, it is necessary to mention some important technical points for which we refer to \cite{ready}.
 \begin{remark}
\label{re:multiply}
Given a piecewise function $\omega$, assumed to be piecewise $C^n$ and positive on $\on$, the set
$\omega\E:=\{\omega F\ |\ F\in \E \}$ is an $(n+1)$-dimensional PW-space on $\on$,  
in which the connection matrices (lower triangular with positive diagonal elements) are given by (see Lemma 39 of \cite{ready})
$$
\mathcal{C}_n(\omega,t_k^+)\ M_k\ \mathcal{C}_{n}(\omega,t_k^-)^{-1}, \quad k=0, \ldots, q,
$$
where, for $x\in [a,b]$ and $\varepsilon \in\{-,+\}$, $\mathcal{C}_{n}(\omega,x^\e)=\bigl(\mathcal{C}_{n}(w,x^\e)_{p,q}\bigr)_{0\leq p,q\leq n}$ stands for the lower triangular square matrix of order $(n+1)$ defined by
$$
\mathcal{C}_{n}(\omega,x^\e)_{p,q}:
=\co{p}{q}\ {w}^{(p-q)}(x^\e),
\quad 0\leq q\leq p\leq n.
$$
\end{remark}
\begin{remark}
\label{re:constants}
The PW-space $\E$ contains constants if and only if, firstly each section-space contains constants, and secondly the first column of each connection matrix is equal to $(1, 0, \ldots, 0)^T$. If so, clearly $\E\subset C^0([a,b])$ and, if $n\geq 1$,  the space $D\E$ is a PW-space on $\on$. The connection matrices in $D\E$ are simply obtained by deleting the first row and column in each $M_k$. 
\end{remark}
\begin{remark}
\label{re:DL0}
Assume that, in the $(n+1)$-dimensional PW-space $\E$, we can find  an element $w_0$ which is positive on $\on$. Then,  denoting by $L_0$ the piecewise division by $w_0$, the space $L_0\E$ contains constants. Combining the previous two reminders shows that, if $n\geq 1$,  the space $DL_0\E$ is an $n$-dimensional PW-space on $\on$.  Therefore, the first order piecewise differential operator $DL_0$ (also named piecewise generalised derivative) diminishes the dimension by  one within the class of all PW-spaces on $\on$. The existence of such a dimension diminishing procedure is thus subject to the existence of a piecewise function $w_0\in\E$ positive on $\on$. This existence is not at all guaranteed, see Section \ref{s:comments}.
\end{remark}

\smallskip
Due to the assumptions on the connection matrices, for any integer $p$, $0\leq p\leq n$, and any $k=1, \ldots, q$, a piecewise function $F$ in the PW-space $\E$ vanishes exactly $p$ times at $t_k^+$ if and only if it vanishes exactly $p$ times at $t_k^-$.  We can therefore count the total number of zeroes of any element of $\E$, including multiplicities up to $(n+1)$. We denote it by $Z_{n+1}(F)$.

\begin{definition}
\label{def:ECP}
The PW-space $\E$ is said to be an {\it Extended Chebyshev Piecewise space} (for short, ECP-space) {\it on $([a,b]; \T)$} when $Z_{n+1}(F)\leq n$ for any non-zero $F\in \E$.
\end{definition}

Equivalently, $\E$ is an ECP-space on $\on$ if, for any positive integers $\mu_1, \ldots, \mu_r$ summing to $(n+1)$, for any pairwise distinct $a_1, \ldots, a_r\in [a,b]$, any convenient $\e_1, \ldots, \e_r\in \{-,+\}$, and any real numbers $\alpha_{i,j}$, $i=1, \ldots, r$, $j=0, \ldots, \mu_{i-1}$, the Hermite interpolation problem
$$
\hbox{find }U\in \E \hbox{ such that  }U^{(j)}(a_i^{\e_i})=\alpha_{i,j} \quad \hbox{for }0\leq j\leq \mu_i-1\hbox{ and for }i=1, \ldots, r,
$$
has a unique solution in $\E$. Clearly, if $\E$ is an ECP-space on $\on$, then, for $k=0, \ldots, q$, any non-zero element $F_k\in\E_k$ satisfies $Z_{n+1}(F_k)\leq n$, which means that the section-space $\E_k$  is an {\it Extended Chebyshev space} (EC-space) {\it on $[t_k, t_{k+1}]$}. Conversely, that each section-space is an EC-space on its own interval does not imply that $\E$ is an ECP-space on $\on$. This justifies the  introduction of the following intermediate definition.

\begin{definition} The PW-space $\E$ is said to be a {\it Piecewise Extended Chebyshev space} (for short, PEC-space) {\it on $([a,b]; \T)$} when for $k=0, \ldots, q$, the section-space $\E_k$ is an EC-space on $[t_k, t_{k+1}]$.
\end{definition}

It is important to mention that, in the situation described in Remark \ref{re:DL0}, a piecewise version of Rolle's theorem says that (see Lemma 38 of \cite{ready})
\begin{equation}
\label{zeroderiv}
Z_n(DL_0F)\geq Z_{n+1}(F)-1 \quad \hbox{for any }F\in\E.
\end{equation}
Accordingly, if $DL_0\E$ is an ($n$-dimensional) ECP-space on $\on$, then $\E$ in turn is an ECP-space on $\on$. Equivalently, the class of all ECP-spaces on $\on$ is closed under continuous integration as well as under multiplication by  (sufficiently piecewise differentiable) positive piecewise functions on $\on$.
A {\it system of piecewise weight functions on $\on$} is a sequence $(w_0,\dots,w_n)$ of piecewise functions on $\on$, such that, for $i=0, \ldots, n$, $w_i$ is positive and $C^{n-i}$ on each $[t_k^+, t_{k+1}^-]$, $k=0, \ldots, q$. With such a system one can associate piecewise generalised derivatives $L_0, \ldots, L_n$ defined in a recursive way as follows:
\begin{equation}
\label{difop}
L_0F:=\frac{F}{w_0}\ ,\quad L_iF:=\frac{1}{w_i}DL_{i-1}F,\quad 1\leq i\leq n.
\end{equation}
By $ECP(w_0, \ldots, w_n)$ we denote the set of all piecewise  functions which are piecewise $C^n$ on $([a,b]; \T)$ and such that $
L_nF$ is constant on $[a,b]$,
with the additional requirement that
$$
L_iF(t_k^+)=L_iF(t_k^-)\hbox{ for }i=0, \ldots, n-1, \hbox{ and for }k=1, \ldots, q.
$$
According to the previous observations, the space $ECP(w_0, \ldots, w_n)$ is an $(n+1)$-dimensional ECP-space on $([a,b];\T)$.
In the special case $q=0$, we recover the well-known procedure to build an $(n+1)$-dimensional EC-space on $[a,b]$ from a
system $(w_0, \ldots, w_n)$ of weight functions on $[a,b]$. The EC-space in question is denoted by $EC(w_0, \ldots, w_n)$.

\smallskip
We conclude this section with two observations.
 \begin{remark}
\label{re:dim1}
Suppose that $n=0$. Then, the one-dimensional PW-space $\E$ on $\on$ is as well a PEC-space on $\on$ or an ECP-space on $\on$. This is clear from the connection conditions (\ref{connection}) and from the fact that being a one-dimensional W-space on $[t_k,t_{k+1}]$  is the same as being a one-dimensional EC-space on $[t_k,t_{k+1}]$. Accordingly, any non-zero $U\in\E$ keeps the same strict sign on $\on$. Equivalently we can state that any one-dimensional PW-space $\E$ on $\on$ can be written as $\E=ECP(w_n)$ where $w_n$ is a piecewise function on $\on$ which is continuous and positive on each $[t_k^+, t_{k+1}^-]$. \end{remark}
\begin{remark}
\label{all}
On account of the previous remark, the class of all spaces of the form $ECP(w_0,  \ldots, w_n)$ coincides with  the class of all PW-spaces on $\on$ in which the dimension diminishing procedure explained in Remark \ref{re:DL0} can be iterated until dimension one.
\end{remark}
%

\subsection{Bernstein-type bases}

We consider again the $(n+1)$-dimensional PW-space $\E$ on $\on$ defined in the previous subsection. Due to the assumptions on the connection matrices, for any $F\in\E$ and any $k=1, \ldots, q$, $F(t_k^+)>0$ if and only  $F(t_k^-)>0$. Along with  the fact that we can count the exact numbers of zeroes at any point in $[a,b]$, this makes the following definitions relevant.

\begin{definition}
\label{BLB}
Given any $c,d\in [a,b]$, $c<d$, and given $B_0, \ldots, B_n\in  \E$, we say that $(B_0, \ldots, B_n)$ is a {\it Bernstein-like basis relative to $(c,d)$} if, for each $i=0, \ldots, n$,  $B_i$ vanishes exactly $i$ times at $c$ and exactly $(n-i)$ times at $d$. We say that it is {\it a positive Bernstein-like basis relative to $(c,d)$} when each $B_i$ additionally satisfies $B_i(x^\e)>0$ for any $x\in]c,d[$.
\end{definition}

\begin{definition}
\label{BB}
A basis  $(B_0, \ldots, B_n)$ in $\E $ is said to be  {\it normalised} if $\sum_{i=0}^nB_i=\un$. A {\it  Bernstein  basis relative to $(c,d)$} is a positive Bernstein-like basis relative to $(c,d)$ which is normalised.
\end{definition}

\medskip
The importance of such bases for ECP-spaces is summarised in the theorem below.



\begin{theorem}
\label{th:ECP_NBLB}
Let $\E$ be an $(n+1)$-dimensional PW-space on $\on$,  assumed to contain constants, with $n\geq 1$. Then, the following properties are equivalent:
\begin{enumerate}[\rm(\roman{enumi})]
\item the space $D\E$ is an ECP-space on $\on$;
\item for any $c,d\in [a,b]$, $c<d$, $D\E$ possesses a Bernstein-like basis relative to $(c,d)$;
\item for any $c,d\in [a,b]$, $c<d$, $\E$ possesses a normalised Bernstein-like basis relative to $(c,d)$;
\item blossoms exist in $\E$. 
\end{enumerate}
\end{theorem}
The equivalence $\rm (i) \Leftrightarrow(ii)$ in Theorem \ref{th:ECP_NBLB} permits to identify all ECP-spaces on $\on$ in the larger class of all PW-spaces on $\on$  by the presence of Bernstein-like bases. For $q=0$, its similarly characterises EC-spaces on $[a,b]$ among all W-spaces on $[a,b]$, see \cite{CA05}. Applied to the section-spaces of $\E$, this yields the following characterisation.

\begin{corollary}
\label{cor:PEC_BLB}
Le $\E$ be an $(n+1)$-dimensional PW-space on $\on$. Then, the following properties are equivalent:
\begin{enumerate}[\rm(\roman{enumi})]
\item the space $\E$ is a PEC-space on $\on$;
\item for any $k=0, \ldots, q$, and any $c,d\in [t_k, t_{k+1}]$, $c<d$, $\E$ possesses a Bernstein-like basis relative to $(c,d)$.
\end{enumerate}
\end{corollary}

The latter theorem and corollary call for a number of important comments below.

\begin{remark}
\label{NBLtoBL}
As is classical, see \cite{CA05}, if $n\geq 1$, and if the PW-space $\E$ possesses a normalised Bernstein-like basis $(B_0, \ldots, B_n)$ relative to $(a,b)$, then this basis generates a Benstein-like basis  relative to $(a,b)$ in the space $D\E$, say $(V_0, \ldots, V_{n-1})$, via the following
formul\ae\
\begin{equation}
\label{eq:Vi}
V_i=D\bigl(B_{i+1}+\cdots+B_n\bigr)=- D\bigl(B_0+\cdots+B_i\bigr), \quad i=0, \ldots, n-1.
\end{equation}
The passage (\ref{eq:Vi}) from normalised Bernstein-like bases to  Bernstein-like bases explains the implication $\rm(iii)\Rightarrow(ii)$ in Theorem \ref{th:ECP_NBLB}.  Moreover, from the expansion $U=\sum_{i=0}^n\alpha_iB_i$ of a function $U\in \E$, we can derive the following expansion of the piecewise function $DU\in D\E$
\begin{equation}
\label{eq:DU}
DU=\sum_{i=0}^{n-1}(\alpha_{i+1}-\alpha_i)V_i.
\end{equation}
\end{remark}
\begin{remark}
\label{re:BLBpositive}
When  (i) of Theorem \ref{th:ECP_NBLB} holds, then all Bernstein-like bases in (ii) of the same theorem can be assumed to be positive in the sense of Definition \ref{BLB}. The same is valid for all Bernstein-like bases in (ii) of Corollary \ref{cor:PEC_BLB}. This is why positivity is often directly included in the definition of Bernstein-like bases. Nevertheless, positivity is not at all needed to prove the implication $\rm(ii)\Rightarrow(i)$ of either Theorem \ref{th:ECP_NBLB} or Corollary \ref{cor:PEC_BLB}, see \cite{CA05}. In the present paper, it  is important to separate the requirement on  the zeroes and the positivity requirement in view of the next section.
\end{remark}
\begin{remark}
\label{re:TP}
Assume that (i) of Theorem \ref{th:ECP_NBLB} holds. Then, according to (iv), blossoms exist in $\E$. Without giving the precise geometrical definition of blossoms, we simply remind the reader that each function $F\in \E^d$ ($d\geq 1$) ``blossoms" into a function $f:[a,b]^n\rightarrow \RR^d$, called the {\it blossom of $F$}.  The three main properties satisfied by blossoms (symmetry, pseudoaffinity in each variable, diagonal property) are the reasons why we can develop a {\it corner cutting} de Casteljau-type algorithm for the evaluation of  the blossom $f$ at any $(x_1, \ldots, x_n)\in [a,b]^n$ as  a convex combination (with coefficients independent of $F$) of the points $f(a^{[n-i]}, b^{[i]})$, $i=0, \ldots, n$, called {\it the B\'ezier points of $F$ relative to $(a,b)$}. The notation $x^{[j]}$ is used with the meaning of $x$ repeated  $j$ times. In particular, this yields:
$$
f(x^{[n]})= F(x)=\sum_{i=0}^nB_i(x)f(a^{[n-i]}, b^{[i]}), \quad \sum_{i=0}^nB_i(x)=1, \quad x\in[a,b].
$$
It follows that $(B_0, \ldots, B_n)$ is the Bernstein basis of $\E$. That it is generated by a corner cutting algorithm guarantees the {\it total positivity} of the Bernstein basis (\ie for any $a\leq x_0<x_1<\cdots <x_n\leq b$, each minor of the  matrix $\bigl(B_i(x_j)\bigr)_{0\leq i,j\leq n}$ is non-negative).
 Concerning the importance of total positivity for geometric design, see \cite{goodman,beyond, CP}. Readers interested in the precise geometrical definition of blossoms can refer to \cite{HP, CA99} for instance.
\end{remark}
Clearly, the property (i) of Theorem \ref{th:ECP_NBLB} implies that $\E$ itself is an ECP-space on $\on$ (note that the converse if not true). This  justifies the definition below.

\begin{definition}
\label{def:good}
A PW-space $\E$ on $\on$ is said to be an {\it ECP-space good for design on $\on$} when it contains constants and when $D\E$ is an ECP-space on $\on$.
\end{definition}

Following from Remark \ref{re:TP} we can state (for details, see \cite{beyond})
\begin{theorem}
Let $\E$ be an $(n+1)$-dimensional ECP-space good for design on $\on$. Then, the normalised Bernstein-like basis of $\E$ relative to $(a,b)$ is its Bernstein basis relative to $(a,b)$ and it is its optimal normalised totally positive basis.
\end{theorem}

\subsection{Dimension diminishing via piecewise generalised derivatives}

Throughout the present section we assume that $n\geq 1$.
In Remark \ref{re:DL0} we have mentioned that diminishing the dimension within the class of all PW-spaces on $\on$ was not always possible. By contrast, this is always possible in the smaller class of all ECP-spaces on $\on$. Below we describe how to proceed. We also examine the same question in the intermediate class of all PEC-spaces on $\on$.

\subsubsection{Within the class of all ECP-spaces on $\on$}

For the proof of the theorem below  readers are referred  to \cite{JAT11, howto}. We simply mention that it is based both on Theorem \ref{th:ECP_NBLB}  and on the properties of blossoms in the ECP-space good for design on $\on$ obtained by continuous integration.

\begin{theorem}
\label{th:w0}
Let $\E$ be an $(n+1)$-dimensional ECP-space on $\on$, and let $(V_0, \ldots, V_n)$ denote a  positive Bernstein-like  basis relative to $(a,b)$. Given a  piecewise function $w_0=\sum_{i=0}^n \alpha_i V_i\in \E$, the following properties are equivalent:
\begin{enumerate}[\rm(\roman{enumi})]
\item $\alpha_0, \ldots,\alpha_n$ are all positive;
\item $w_0$ is positive on $\on$, and if $L_0$ denotes the piecewise division by $w_0$, the $(n+1)$-dimensional space $L_0\E$  is an ECP-space good for design on $\on$.
\end{enumerate}
\end{theorem}

Property (ii) of Theorem \ref{th:w0} means that the $n$-dimensional PW-space $DL_0\E$ is an ECP-space on $\on$. In other words, when (i) holds,  the first order piecewise differential operator $DL_0$ diminishes the dimension by one within the class of all ECP-spaces on $\on$. Iterating the process leads to

\begin{corollary}
\label{cor:allECP}
Any  $(n+1)$-dimensional ECP-space on $\on$ is of the form $ECP(w_0,  \ldots, w_n)$, where $(w_0, \ldots, w_n)$  is a system of piecewise weight functions on $\on$.
\end{corollary}
On account of Definition \ref{def:good}, Corollary \ref{cor:allECP} obviously implies that any $(n+1)$-dimensional ECP-space good for design on $\on$ is of the form $ECP(\un, w_1, \ldots, w_n)$.

\begin{remark}
\label{manyL1}
As a consequence of Theorem \ref{th:w0}, as soon as $n\geq 1$, the dimension diminishing process produces infinitely many different $n$-dimensional  ECP-spaces $DL_0\E$.
\end{remark}

\begin{remark}
\label{re:BLtoBL}
Assuming that (i) of Theorem \ref{th:w0} holds, clearly the Bernstein basis $(B_0, \ldots, B_n)$ relative to $(a,b)$ in the space $L_0\E$ is given by
$$
B_i:=\frac{ \alpha_i V_i}{w_0}, \quad i=0, \ldots, n.
$$
This basis generates in turn a positive  Bernstein-like  basis $(\overline V_0, \ldots, \overline V_{n-1})$ relative to $(a,b)$ in the ECP-space $DL_0\E$ via (\ref{eq:Vi}). This recursive construction of bases will be essential in the next section.
\end{remark}

\subsubsection{Within the class of all PEC-spaces on $\on$}
\label{L1PEC}

Observe that, if $q=0$, Theorem \ref{th:w0} describes how to diminish the dimension in the class of all EC-spaces $\E$ on $[a,b]$, that  is, all possible first steps to write $\E$ as $\E=EC(w_0,  \ldots, w_n)$. We would like to draw the reader's attention that this is not valid for EC-spaces on intervals which are not assumed to be closed and bounded.

Here we assume  that $\E$ is a PEC-space on $\on$.
How can we similarly diminish the dimension within the class of all PEC-spaces on $\on$? This is recalled in Proposition below, as  a straightforward consequence of the non-piecewise version of Theorem \ref{th:w0}.

Let us first observe that, a priori, we do not have a global (positive) Bernstein-like basis (that is, relative to $(a,b)$) in the space $\E$. By contrast,  as recalled in Corollary \ref{cor:PEC_BLB} and Remark \ref{re:BLBpositive} we know that, for each $k=0, \ldots, q$, the space  $\E$ possesses a positive Bernstein-like basis relative to $(t_k, t_{k+1})$. We denote such a basis by  $(V_{k,0}, V_{k,1}, \ldots, V_{k,n})$ and we refer to it as {\it a $k$th local positive Bernstein-like basis} in $\E$.

\begin{proposition}
\label{prop:DL0PEC}
Assume that $\E$ is an $(n+1)$-dimensional PEC-space on $\on$. Consider a piecewise function $w_0\in \E$, expanded in the local positive Bernstein-like bases of $\E$ as
\begin{equation}
\label{w0local}
w_0=\sum_{r=0}^n\delta_{k,r}V_{k,r}, \quad k=0, \ldots, q.
\end{equation}
Then, the following two properties are equivalent
\begin{enumerate}[\rm(\roman{enumi})]
\item all coefficients $\delta_{k,r}$, $r=0, \ldots, n$, $k=0, \ldots, q$, are positive;
\item  $w_0$ is positive on $\on$ and if $L_0$ denotes the division by $w_0$, $DL_0\E$ is an ($n$-dimensional) PEC-space on $\on$.
\end{enumerate}
\end{proposition}

It is important to emphasise that, in the PEC context, the existence of a function $w_0$ satisfying (i) is not guaranteed, see Subsection \ref{s:comments}.
\begin{remark}
\label{re:localpassage}
Assuming that (i) of Proposition \ref{prop:DL0PEC} holds, the passage from local positive Bernstein-like bases in $\E$ to local Bernstein bases in $L_0\E$, and then  to local positive Bernstein-like bases in $DL_0\E$ can be described by local versions of Remark \ref{re:BLtoBL}.
\end{remark}

\subsubsection{Concluding remarks}
\label{s:comments}
The three classes of PW-, PEC-, ECP-spaces on $\on$ are closed under continuous integration on $[a,b]$ and under multiplication by positive piecewise functions on $\on$  (assumed to be sufficiently piecewisely differentiable).

Only in the largest class of PW-spaces on $\on$ diminishing by one the dimension via a piecewise generalised derivative  is equivalent to finding a positive piecewise function $w_0$ in  the considered space. In each of the two other classes the positivity requirement is not sufficient in general to remain in the same class. Simple counterexamples can be found in \cite{JAT11} for instance.

Only in the smallest class of ECP-spaces on $\on$ the dimension diminishing process via piecewise generalised derivatives is always possible. This can be easily illustrated by the following classical example. Given $h>0$, take $t_0:=-h$, $t_1:=0$, $t_2:=h$, and $[a,b]:=[t_0, t_2]$. Let $\E$ be the space spanned over $[a,b]$ by the two functions $\cos, \sin$. It can be considered as well as a PW-space on $\on$, with $\T:=(t_1)$, with section-spaces its restrictions $\E_0, \E_1$ to $[t_0, t_1]$ and $[t_1, t_2]$, respectively, and with the identity connection matrix of order two at $t_1$. As is well known, each section-space is an EC-space on its interval if and only if $h<\pi$. This is therefore the necessary and sufficient condition for $\E$ to be a PEC-space on $\on$. Now, $\E$ can also be described as the set of all functions on $[a,b]$ which are of the form $\alpha \cos(x-\beta)$, for some $\alpha, \beta\in \RR$. Accordingly, a non-zero $w_0$ exists in $\E$ if and only if $h<\frac{\pi}{2}$, that is, if and only if $\E$ is an EC on $[a,b]$, or equivalently, if and only if $\E$  an ECP-space on $\on$. Note that, here, when $2h<\pi$, the  requirement (i) of Theorem \ref{th:w0} is equivalent to $w_0$ being positive because the space is two-dimensional.

\section{How to test if a given PW-space is an ECP-space}
\label{s:test}
Among all PW-spaces on $\on$ containing constants, we would like to be able to recognize those which are ECP-space good for design on $\on$. On account of Definition \ref{def:good} this amounts to being able to answer the following question:  an $(n+1)$-dimensional PW-space $\E$ on $\on$ being given:
\begin{equation}
\label{question}
\hbox{is }\E \hbox{ an ECP-space on }\on?
\end{equation}
Trying to answer Question (\ref{question}) is the object of the present section.
\subsection{The guiding principles}

Clearly, on account of the $(q+1)$ different section-spaces and of the presence of connection matrices, as soon as $q\geq 1$, this question is difficult to answer, even in rather low dimensions. Still, in Remark \ref{all} we have given an interesting  general theoretical answer to (\ref{question}):  the given PW-space $\E$  is an ECP-space on $\on$ if and only if we are able to iteratively build first order piecewise differential operators so as to diminish the dimension by one at each step, until we reach dimension one. This will be our guiding principle.

\smallskip
The problem is that, in the large class of PW-spaces, we have no constructive way to test if it is possible to find such piecewise differential operators. 
Now, a first necessary condition to give a positive answer to (\ref{question}) is that each section-space should be an EC-space on its own interval. How to test this point is another problem for which we refer to \cite{CMP, JJA12, CMP1} for instance. Therefore, from now on we will directly assume that {\bf $\E$ is a given  $(n+1)$-dimensional PEC-space on $\on$}.  The interest of working in this framework is that we will have at our disposal local positive Bernstein-like bases, see Subsection \ref{L1PEC}.  Proposition \ref{prop:DL0PEC} shows that these local bases are practical  tools to test whether or not a given function $w_0\in \E$ enables us to diminish the dimension within the class of all PEC-spaces on $\on$.

Now, how to select the functions $w_0\in \E$ which  should be tested? As a matter of fact, the possible candidates will be selected by analogy with the case where we do know that $\E$ is an ECP-space on $\on$, that is, by analogy with (i) of Theorem
\ref{th:w0}. Nonetheless, in order to imitate condition (i) of Theorem
\ref{th:w0}, we first need to have at our disposal a basis $(V_0, \ldots, V_n)$ in $\E$ which resembles a positive Bernstein-like relative to $(a,b)$. This point will thus be a preliminary requirement before trying to diminish the dimension.

Finally, though (i) of Theorem \ref{th:w0} offers infinitely many different possibilities to select a function $w_0$, it will actually be sufficient to test one single appropriately chosen function $w_0\in \E$. The various steps of the test will be based on specific properties in ECP-spaces on $\on$ that we present in the next subsection.
\subsection{Useful necessary conditions}
Throughout the present subsection we assume that $\E$ is an ECP-space on $\on$. Let $(V_0, \ldots, V_n)$ be a given positive Bernstein-like basis relative to $(a,b)$ in the space $\E$. It will be convenient to refer to it as {\it the initial  global positive Bernstein-like basis}.
\subsubsection{Some properties of the global Bernstein-like basis}

Below we point out a few properties  of the basis $(V_0, \ldots, V_n)$ which directly follow from the definition of positive Bernstein-like bases.  We state them as lemmas because they will serve to build the test later on.

\smallskip
For each $i=0, \ldots, n$, we know that the function $V_i$ is positive on $]a,b[$ and that it vanishes $i$ times at $a$ and $(n-i)$ times at $b$. Its Taylor expansions at $a$ and $b$ show that
\begin{lemma}
\label{lem:positivity_ab}
The global positive Bernstein-like  basis $(V_0, \ldots, V_n)$ satisfies
\begin{equation}
\label{positivity_ab}
{V_i}^{(i)}(a)>0 \quad\hbox{and}\quad (-1)^{n-i}\,{V_i}^{(n-i)}(b)>0, \quad i=0, \ldots, n.
\end{equation}
\end{lemma}
\noindent
We refer to (\ref{positivity_ab}) as the {\it endpoint positivity conditions}.

\smallskip
For $k=0, \ldots, q$, let $(V_{k,0}, \ldots, V_{k,n})$ denote a given $k$th local positive Bernstein-like basis in $\E$. We can then consider the  expansions of  the global basis $(V_0, \ldots, V_n)$ in the local ones, which we denote as
\begin{equation}
\label{Vi_local}
V_i=\sum_{r=0}^n\gamma_{i,k,r}V_{k,r}, \quad k=0, \ldots, q, \quad i=0, \ldots, n.
\end{equation}
Without loss of generality we can assume that the global basis $(V_0, \ldots, V_n)$ is derived from the global Bernstein basis in the $(n+2)$-dimensional ECP-space $\widehat \E$ good for design on $\on$ obtained by continuous integration of $\E$,  via the classical procedure (\ref{eq:Vi}). The de Casteljau evaluation algorithm for blossoms in $\widehat\E$ and the local formul\ae\ (\ref{eq:DU}) yield the following result (for details see Proposition 6.8 of \cite{CAGD2016}):
%
\begin{lemma}
\label{lem:pos_local_coef}
The coefficients of the expansions (\ref{Vi_local}) satisfy
\begin{equation}
\label{eq:pos_local_coef}
	\begin{array}{lll}
	&\gamma_{i,0,0}=\gamma_{i,0,1}= \cdots =\gamma_{i,0,i-1}=0; \quad \gamma_{i,0,r}>0 \quad\hbox{for }r=i, \ldots, n,\\
	\\
	&\gamma_{i,k,r}>0  \quad\hbox{for }r=0, \ldots, n, \quad\hbox{when }1\leq k\leq q-1,\\
	\\
	&\gamma_{i,q,r}>0 \quad\hbox{for }r=0, \ldots, i, \quad \gamma_{i,q,i+1}=\cdots=\gamma_{i,q,n-1}=\gamma_{i,q,n}=0.
	\end{array}
 \end{equation}
\end{lemma}
%

%
\subsubsection{One specific construction of piecewise generalised derivatives}
\label{onesystem}%

In Remark \ref{manyL1} we have reminded the reader that, if $n\geq 1$, the space $\E$ generates  infinitely many different $n$-dimensional ECP-spaces $DL_0\E$ on $\on$ via the dimension diminishing process described in Theorem \ref{th:w0}.  In this subsection we will describe one specific systematic way to build a system $(w_0, \ldots, w_n)$ of piecewise weight functions on $\on$ such that $\E=ECP(w_0, w_1, \ldots, w_n)$. According to Theorem \ref{th:w0}, since  the function
\begin{equation}
\label{eq:w0}
w_0:=V_0+V_1+\cdots+V_n.
 \end{equation}
has clearly positive coordinates in the global positive Bernstein-like basis, it can be taken as a first piecewise weight function on $\on$. If $L_0$ denotes the piecewise division by $w_0$, we know that the global  Bernstein basis in $L_0\E$  is given by
\begin{equation}
\label{eq:Bi}
B_i=\frac{V_i}{w_0}, \quad i=0, \ldots,n.
 \end{equation}
Moreover, according to (\ref{eq:Vi}),  the sequence $(\overline V_0, \overline V_1, \ldots, \overline V_{n-1})$ defined by
\begin{equation}
\label{eq:Vibar}
\overline V_i:=D(B_{i+1}+\cdots+B_n), \quad i=0, \ldots, n-1,
 \end{equation}
 is a global positive Bernstein-like basis in the ECP-space $DL_0\E$ on $\on$. Then,
\begin{equation}
\label{eq:w1}
w_1:=\overline V_0+ \overline V_1+\cdots+ \overline V_{n-1}
 \end{equation}
is taken as the next piecewise weight function, and with $L_1$ defined as in (\ref{difop}), the $n$-dimensional space $L_1\E$ is an ECP-space good for design on $\on$ and so forth. This describes one simple systematic procedure to write  $\E$  as $\E=ECP(w_0, \ldots, w_n)$.


%
\subsubsection{Iterated local expansions}
\label{iterated}
In the previous iterative construction of the system $(w_0, \ldots, w_n)$ of piecewise weight functions on $\on$, all global bases are derived form the initial global positive Bernstein-like basis. The procedure induces a simultaneous iterated construction of local bases from the initial  local positive Bernstein-like bases. As a consequence, the local expansions of the successive global positive Bernstein-like bases can be derived from the initial one (\ref{Vi_local}). One step of this iteration is subsequently described in detail.

\smallskip
Take $k\in\{0, \ldots, q\}$. On account of (\ref{Vi_local}), the $k$th local expansion of the piecewise weight function $w_0$ introduced in (\ref{eq:w0}) is
\begin{equation}
\label{eq:w0_local}
w_0=\sum_{r=0}^n\delta_{k,r}V_{k,r}, \quad\hbox{with }\delta_{k,r}:=\sum_{i=0}^n\gamma_{i,k,r}.
 \end{equation}
In $L_0\E$, the $k$th local Bernstein basis $(B_{k,0}, \ldots, B_{k,n})$ is thus given by
\begin{equation}
\label{eq:BB_L0}
B_{k,r}:=\frac{\delta_{k,r}V_{k,r}}{w_0}, \quad r=0, \ldots, n.
 \end{equation}
The functions
\begin{equation}
\label{eq:localBLBbar}
\overline V_{k,r}:=D(B_{k,r+1}+\ldots+ B_{k,n}), \quad r=0, \ldots, n-1,
 \end{equation}
form a $k$th local positive Bernstein-like basis in $DL_0\E$. This local basis will be used to obtain the $k$th local expansion
\begin{equation}
\label{eq:Vibar_local}
\overline V_i=
\sum_{r=0}^{n-1}\overline\gamma_{i,k,r}\overline V_{k,r}, \quad i=0, \ldots, n-1,
 \end{equation}
of the global positive Bernstein-like basis $(\overline V_0, \ldots, \overline V_{n-1})$ of $DL_0\E$ obtained in (\ref{eq:Vibar}).  We know that this expansion will be derived from the expansion of the function $B_{i+1}+\ldots+ B_{n}$ in the local Bernstein basis $(B_{k,0}, \ldots, B_{k,n})$ via formula (\ref{eq:DU}). Now, from (\ref{eq:Bi}), (\ref{Vi_local}), and (\ref{eq:BB_L0}) we can write:
\begin{equation}
\label{eq:sum_B}
B_{i+1}+\ldots+ B_{n}=\frac{1}{w_0}\left(V_{i+1}+\cdots+V_n\right)=\frac{1}{w_0}\left(\sum_{j=i+1}^n\gamma_{j,k,r}\frac{w_{0}B_{k,r}}{\delta_{k,r}}\right).
 \end{equation}
Accordingly, (\ref{eq:DU}) yields
\begin{equation}
\label{eq:gamma_bar}
\overline\gamma_{i,k,r}=\frac{\sum_{j=i+1}^n\gamma_{j,k,r+1}}{\delta_{k,r+1}}-\frac{\sum_{j=i+1}^n\gamma_{j,k,r}}{\delta_{k,r}}, \quad 0\leq i,r\leq n-1.
 \end{equation}

\subsection{The numerical test}

Throughout this section we assume that $\E$ is an $(n+1)$-dimensional PEC-space on $\on$. For $k=0, \ldots, q$, we denote by $(V_{k,0}, \ldots, V_{k,n})$ a given $k$th local positive Bernstein-like basis  in $\E$ (see Corollary \ref{cor:PEC_BLB} and Remark \ref{re:BLBpositive}).

\subsubsection{Step 0}
\label{test1}
For $\E$ to be an ECP-space on $\on$ it is necessary that $\E$ possess a global Bernstein-like basis satisfying the endpoint positivity conditions (\ref{positivity_ab}) and whose local expansions satisfy (\ref{eq:pos_local_coef}).

\smallskip
\noindent $\bullet $ {\bf Part 1-}
The first part of the test is:
{\it \begin{itemize}
\item[1)] Given $i=0, \ldots, n$, does  there exist a piecewise function $V_i\in\E$ satisfying the $n$  Hermite interpolating conditions
\begin{equation}
\label{Hermite_endcond1}
V_i(a)=V'_i(a)=\dots=V_i^{(i-1)}(a)=0, \qquad V_i(b)=V'_i(b)=\dots=V_i^{(n-i-1)}(b)=0,
\end{equation}
along with
\begin{equation}
\label{Hermite_endcond2}
V_i^{(i)}(a)=1 \quad \mathrm{if} \quad i \leq \left\lfloor\frac{n}{2}\right\rfloor, \qquad V_i^{(n-i)}(b)=(-1)^{n-i} \quad \mathrm{if} \quad i > \left\lfloor\frac{n}{2}\right\rfloor?
\end{equation}
\end{itemize}
}
To this end, for each $i=0,\dots,n$, we set up the linear system of size $(q+1)(n+1)$ satisfied by the  coefficients $\gamma_{i,k,r}$, $k=0, \ldots, q$, $r=0, \ldots, n$, of the local expansions (\ref{Vi_local}) of the piecewise $V_i$ we are looking for. The $q(n+1)$ equations expressing the connection conditions (\ref{connection}) are completed by the $(n+1)$ conditions (\ref{Hermite_endcond1}) and (\ref{Hermite_endcond2}).  We should observe that the $n$ equations due to  (\ref{Hermite_endcond1}) are as follows:
\begin{equation}
\label{eq:pos_local_end1}
	\gamma_{i,0, 0}=\gamma_{i,0,1}= \cdots =\gamma_{i,0,i-1}=0,\quad
	\gamma_{i,q,i+1}=\cdots=\gamma_{i,q,n-1}=\gamma_{i,q,n}=0, 
 \end{equation}
and  the last one due to  (\ref{Hermite_endcond2}) can be written as:
\begin{equation}
\label{eq:pos_local_end2}
\gamma_{i,0,i}=\frac{1}{{V_{0,i}}^{(i)}(a)}\quad\hbox{if } i \leq \left\lfloor\frac{n}{2}\right\rfloor,  \qquad \gamma_{i,q,i}=\frac{(-1)^{n-i}}{{V_{q,i}}^{(n-i)}(b)}\quad\hbox{if } i > \left\lfloor\frac{n}{2}\right\rfloor.
 \end{equation}

In case there exists an integer $i$ such that the system in the unknowns $\gamma_{i,k,r}$, $k=0, \ldots, q$, $r=0, \ldots, n$, has no solution,  we can definitively say that $\E$ is not an ECP-space on $\on$.

Assume that the answer to Question 1) is positive for each $i=0, \ldots, n$. Then the sequence $(V_0, \ldots, V_n)$ is a basis of $\E$, and it satisfies half of the endpoint positivity conditions  (\ref{positivity_ab}). It is therefore a good candidate to be a global positive Bernstein-like basis but we cannot yet assert that it is indeed such a basis.

\medskip
\noindent $\bullet $ {\bf Part 2: }
\label{test2}
What we test now is:

{\it \begin{itemize}
\item[2)] Do all sequences $(\gamma_{i,k,0}, \ldots, \gamma_{i,k, n})$, $i=0, \ldots,n$, $k=0, \ldots, q$, satisfy  (\ref{eq:pos_local_coef})?
\end{itemize}
}
\noindent
A negative answer to Question 2) means that $\E$ is not an ECP-space on $\on$. Assume that we obtain a positive answer. Then, the local expansions prove that $(V_0, \ldots, V_n)$ is a global positive Bernstein-like basis in $\E$.
Moreover, we can assert that
$$
\delta_{k,r}:=\sum_{i=0}^n\gamma_{i,k, r}>0 \quad \hbox{for each }r=0, \ldots, n, \quad \hbox{and for each }k=0, \ldots, q.
$$
Accordingly, the positive piecewise function
$$
w_0:=V_0+V_1+\dots+V_n\in \E,
$$
satisfies (i) of Proposition \ref{prop:DL0PEC} since it can be expanded as in (\ref{eq:w0_local}).
The corresponding space $DL_0\E$  is thus an $n$-dimensional PEC-space on $\on$.

\subsubsection{Iteration of the test}
\label{iter}
Let us assume that  the two questions 1) and 2) above have received positive answers. Then, we can iterate the test, that is, we can apply it to the
$n$-dimensional PEC-space $DL_0\E$ on $\on$.

Before starting the iteration, it is necessary to give further comments on the passage from $\E$ to $L_0\E$, for which we can copy and paste every single formula of Sections \ref{onesystem} and \ref{iterated}. In particular, the basis $(B_0, \ldots, B_n)$ defined in (\ref{eq:Bi})  is a global Bernstein basis in the PEC-space $L_0\E$. The section-spaces of $L_0\E$ being EC-spaces good for design on their intervals, for each $k=0, \ldots, q$, the basis $(B_{k,0}, \ldots, B_{k,n})$ introduced in (\ref{eq:BB_L0}) is its $k$th local Bernstein basis. Similarly, the functions $\overline V_{k,r}$, $r=0, \ldots, n-1$, defined by  (\ref{eq:localBLBbar}) form a $k$th local positive Bernstein-like basis in the PEC-space $DL_0\E$.

By contrast, it should be observed that, in the PEC-space $DL_0\E$, the basis $\left(\overline V_0, \overline V_1, \ldots, \overline V_{n-1}\right)$ defined in (\ref{eq:Vibar})  is not necessarily a global positive Bernstein-like  basis. Here, we can only assert that it is a global Bernstein-like basis. Nevertheless, it automatically satisfies the adapted end point positivity conditions (\ref{positivity_ab}) due to the positivity of the piecewise function $w_0$.

Accordingly,  the first part of the test is useless. We only have to apply the second part, thus applying Question 2) to all coefficients $\overline\gamma_{i,k,r}$ in the expansions  (\ref{eq:Vibar_local}). As in Section \ref{iterated}, we know that these coefficients are given by (\ref{eq:gamma_bar}), the $\delta$'s being taken from (\ref{eq:w0_local}). If the answer is negative, we stop the test: $\E$ is not an ECP-space on $\on$. Otherwise we continue the dimension diminishing procedure.

\smallskip
Similarly, if, at some further iteration, Question 2) is answered negatively, the test stops: $\E$ is not an ECP-space on $\on$.
Suppose that we can perform $n$ such dimension diminishing steps with positive answers to  Question 2)  at each step,  via positive piecewise functions $w_0$, $w_1, \ldots, w_{n-1}$ on $\on$, successively defined in accordance with (\ref{eq:w0}) and (\ref{eq:w1}). After these $n$ steps, we obtain a one-dimensional PEC-space on $\on$, $DL_{n-1}\E$ -- where $L_{n-1}$ is defined as in (\ref{difop}) -- possessing a global Bernstein-like basis, say $\mathcal{V}_0$, satisfying  the endpoint positivity conditions $\mathcal{V}_0(a)>0$. $\mathcal{V}_0(b)>0$. From Remark \ref{re:dim1} we know that $DL_{n-1}\E=ECP(w_n)$ where $w_n:=\mathcal{V}_0$ is our last piecewise weight function.
Then, we can state that $\E=ECP(w_0, \ldots, w_n)$. In retrospect we can say that each global Bernstein-like basis  in the process is indeed a global positive  Bernstein-like basis. Accordingly, the test  is successful as soon as we reach dimension two with a positive answer to Question 2).

\subsubsection{Implementation}

Following the theoretical description given in Subsections 3.3.1 and 3.3.2,  the implementation of the numerical test comprises the following steps.

\medskip
\noindent $\bullet $ {\bf T0:}
This part of the test corresponds to Step 0.
\begin{enumerate}
\item[T0.1]
Given $i=0, \ldots, n$,
if the square  matrix of order $(q+1)(n+1)$ involved in the computation of the coefficients $\gamma_{i,k,r}$, $k=0, \ldots, q$, $r=0, \ldots, n$, is nearly-singular or/and very ill-conditioned, we stop the test: in that case the system is not solvable or the solution cannot be computed with sufficient accuracy. Otherwise we solve the system.  If all the functions $V_i$, $i=0,\dots,n$ can be computed, we proceed to T0.2.
\item[T0.2] Test if all $\gamma_{i,k,r}$ other than those in (\ref{eq:pos_local_end1}) and (\ref{eq:pos_local_end2}) are positive. If the test fails, stop. Otherwise, with $p=0$  and $\gamma_{i,k,r}^{\{0\}}\coloneqq \gamma_{i,k,r}$ for $k=0, \ldots, q$, $i,r=0, \ldots, n$ compute (see (\ref{eq:gamma_bar}) and (\ref{eq:w0_local}))
\begin{equation}
\label{eq:gamma_p}
\gamma_{i,k,r}^{\{p+1\}}=\frac{\sum_{j=i+1}^{n-p}\gamma_{j,k,r+1}^{\{p\}}}{\sum_{j=0}^{n-p}\gamma_{j,k,r+1}^{\{p\}}}-\frac{\sum_{j=i+1}^{n-p}\gamma_{j,k,r}^{\{p\}}}{\sum_{j=0}^{n-p}\gamma_{j,k,r}^{\{p\}}}, \quad 0\leq i,r\leq n-p-1, \ k=0, \ldots, q,
 \end{equation}
 and proceed to T1.
\end{enumerate}
\medskip
\noindent $\bullet $ {\bf T1:}
This part of the test corresponds to the dimension diminishing process. To describe it, for some integer $p$, $1\leq p\leq n-1$, consider real numbers $\gamma_{i,k,r}^{\{p\}}$, $i,r=0,\dots,n-p$, $k=0,\dots,q$, such that
$$\gamma_{i,0,0}^{\{p\}}=\cdots= \gamma_{i,0,i-1}^{\{p\}}=0<\gamma_{i,0,i}^{\{p\}}, \qquad
\gamma_{i,q,i}^{\{p\}}>0= \gamma_{i,q,i+1}^{\{p\}}=\cdots=\gamma_{i,q,n-p}^{\{p\}}.
$$
\begin{enumerate}
\item[T1.1] Test if all other $\gamma_{i,k,r}^{\{p\}}$ are positive.  As soon as the positivity test fails, stop.
\item[T1.2] If T1.1 is successful, compute $\gamma_{i,k,r}^{\{p+1\}}$, $i,r=0,\dots,n-p-1$, $k=0,\dots,q$, according to formula \eqref{eq:gamma_p}.
Return to T1.1 with $p$ replaced by $(p+1)$. 
\end{enumerate}
The test T1 starts at $p=1$. For numerical reasons, all positivity tests are actually replaced by $\gamma_{i,k,r}^{\{p\}}>$ tol, where tol is a small positive number (in our experiments, tol $= 1e-30$).
When the test stops  either at T0 or at T1 for some $p\leq n-1$, then, numerically speaking, we consider that the space $\E$ is not an ECP-space on $\on$.  When the test successfully continues to T1 until $p=n-1$, numerically speaking, we consider that the space $\E$ is an ECP-space on $\on$.

\medskip

The following MATLAB function takes as input a matrix \texttt{gamma} of dimension $(n+1)\times (q+1)\times(n+1)$, where \texttt{gamma(i,k,r)}$=\gamma_{i,k,r}$ are the coefficients in equation \eqref{Vi_local}.
It returns a variable \texttt{test}, which is equal to zero if the test on the positivity of the local coefficients fails at some step.
The loops in the variables \texttt{k}, \texttt{i}, \texttt{r} iterate respectively over the intervals, the
elements of the global bases and the coefficients of their local expansions on each interval.

{\small
\begin{Verbatim}[commandchars=\\\{\}]
function test = ECP_test(gamma,tol)
n=size(gamma,1)-1;
q=size(gamma,2);
test=1;
p=0;
while (p<=n-1 & test)
    np=n+1-p;
    % loop over all intervals to test whether the local coefficients gamma(i,k,r)
    % for the (n+1-j)-dimensional PEC-space satisfy (10)
    for k=1:q
        for i=1:np
            for r=1+(k==1)*(i-1):np-(k==q)*(np-i)
                if (gamma(i,k,r) < tol)
                    test=0;
                    return
                end
            end
        end
    end
    % loop over all intervals to compute the local coefficients gamma(i,k,r)
    % for the the (n-p)-dimensional PEC-space through (20)
    for k=1:q
        % sums in (20)
        for i=np-1:-1:1
            for r=1:np
                gamma(i,k,r)=gamma(i,k,r)+gamma(i+1,k,r);
            end
        end
        % divisions in (20)
        for r=1:np
            gamma1kr=gamma(1,k,r);
            for i=2:np
                gamma(i-1,k,r)=gamma(i,k,r)/gamma1kr;
            end
        end
        % differences in (20)
        for i=1:np-1
            for r=1:np-1
                gamma(i,k,r)=gamma(i,k,r+1)-gamma(i,k,r);
            end
        end
    end
    p=p+1;
end
\end{Verbatim}
}

\begin{figure}[h]
\centering
\includegraphics[width=0.55\textwidth,trim = 0 0cm 0 0cm, clip]{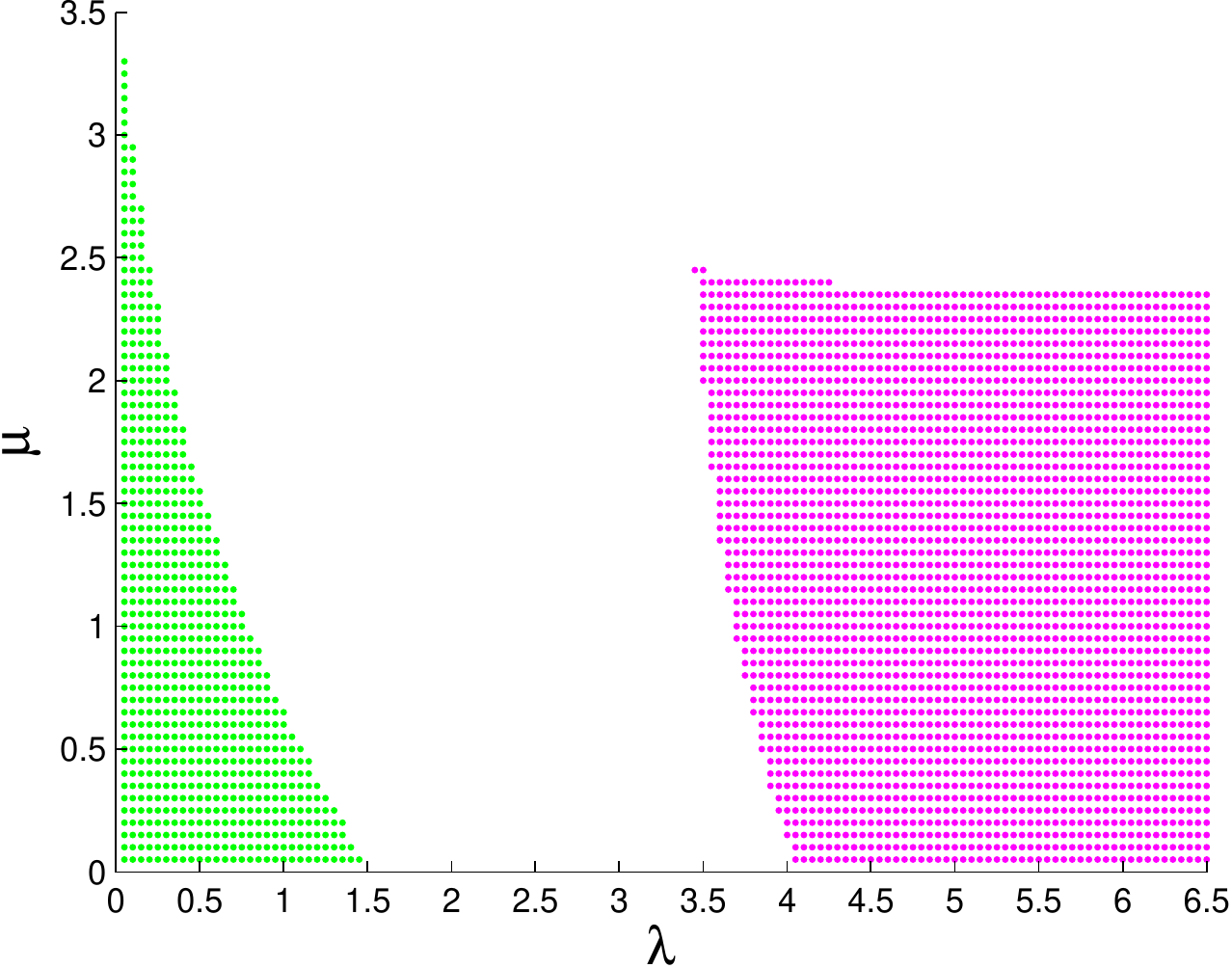}
\caption{Graphical illustration of the output of the numerical test for the THTH PEC-space in Section \ref{THTH}.}
\label{fig:gfd_TH}
\end{figure}
\begin{table}[h]
\centering
\resizebox{\textwidth}{!}{
\begin{tabular}{| c | c c |}
\cline{2-3}
\multicolumn{1}{c|}{} & $\E$  & $DL_0\E$ \\
\hline
\raisebox{-0.25cm}{\rotatebox{90}{k=0}} &
$
\begin{array}{ccccc}
{} & {} & {} & {} & {}\\[0.3ex]
1 & 3.3817 & 8.9847 & 2.6979 & 0.62437 \\[0.3ex]
0 & 1 & 4.8671 & 2.0399 & 0.59032 \\[0.3ex]
0 & 0 & 1 & 0.84793 & 0.37871 \\[0.3ex]
0 & 0 & 0 & 0.0035679 & 0.0036312 \\[0.3ex]
0 & 0 & 0 & 0 & 2.8565e-06 \\[3ex]
\end{array}
$
&
$
\begin{array}{cccc}
0.22822 & 0.16682 & 0.12226 & 0.091739 \\[0.3ex]
0 & 0.067332 & 0.085013 & 0.087064 \\[0.3ex]
0 & 0 & 0.00063834 & 0.0016372 \\[0.3ex]
0 & 0 & 0 & 1.7886e-06 \\[0.3ex]
\end{array}
$
\\
\hline
\raisebox{-0.25cm}{\rotatebox{90}{k=1}} &
$
\begin{array}{ccccc}
{} & {} & {} & {} & {}\\[0.3ex]
0.62437 & 0.46973 & 0.20328 & 0.24318 & 0.16292 \\[0.3ex]
0.59032 & 0.95505 & 0.41134 & 0.60222 & 0.48348 \\[0.3ex]
0.37871 & 1.0734 & 1.6102 & 3.4078 & 3.1382 \\[0.3ex]
0.0036312 & 0.014855 & 0.040202 & 0.33095 & 0.41968 \\[0.3ex]
2.8565e-06 & 1.4492e-05 & 5.2756e-05 & 0.00069584 & 0.0069185 \\[3ex]
\end{array}
$
&
\hspace{0.7cm}
$
\begin{array}{cccc}
0.20404 & 0.097174 & 0.036704 & 0.014351 \\[0.3ex]
0.19364 & 0.29561 & 0.086956 & 0.030893 \\[0.3ex]
0.0036412 & 0.011855 & 0.054563 & 0.028967 \\[0.3ex]
3.9781e-06 & 1.7524e-05 & 0.00012848 & 0.0014911 \\[0.3ex]
\end{array}
$
\\
\hline
\raisebox{-0.25cm}{\rotatebox{90}{k=2}} &
$
\begin{array}{ccccc}
{} & {} & {} & {} & {}\\[0.3ex]
0.16292 & 0.024924 & 0.0078407 & 0.0061688 & 0.0037487 \\[0.3ex]
0.48348 & 0.19339 & 0.094695 & 0.16522 & 0.17542 \\[0.3ex]
3.1382 & 1.7563 & 0.65982 & 1.7289 & 2.7985 \\[0.3ex]
0.41968 & 0.35967 & 0.18827 & 0.37974 & 0.7234 \\[0.3ex]
0.0069185 & 0.010689 & 0.012989 & 0.040777 & 0.26258 \\[3ex]
\end{array}
$
&
$
\begin{array}{cccc}
0.028059 & 0.0024918 & 0.0054787 & 0.0017123 \\[0.3ex]
0.060401 & -0.013313 & 0.032559 & 0.028646 \\[0.3ex]
0.056635 & 0.050919 & -0.027658 & 0.06756 \\[0.3ex]
0.0029153 & 0.0089212 & 0.0040909 & 0.048678 \\[0.3ex]
\end{array}
$
\\
\hline
\raisebox{-0.25cm}{\rotatebox{90}{k=3}} &
$
\begin{array}{ccccc}
{} & {} & {} & {} & {}\\[0.3ex]
   0.0037487 & 0 & 0 & 0 & 0 \\[0.3ex]
0.17542 & 0.12345 & 0 & 0 & 0 \\[0.3ex]
2.7985 & 2.8763 & 1.1487 & 0 & 0 \\[0.3ex]
0.7234 & 0.81068 & 0.46275 & 1 & 0 \\[0.3ex]
0.26258 & 0.39133 & 0.32664 & 0.75473 & 1 \\[3ex]
\end{array}
$
&
$
\begin{array}{cccc}
 0.00094577 & 0 & 0 & 0 \\[0.3ex]
0.015822 & 0.029382 & 0 & 0 \\[0.3ex]
0.037316 & 0.12123 & 0.5927 & 0 \\[0.3ex]
0.026886 & 0.075403 & 0.26157 & 0.56989 \\[0.3ex]
\end{array}
$
\\
\hline
\end{tabular}
}
\caption{Coefficients of the local expansion of the Bernstein-like bases for the THTH PEC-space  in Section \ref{THTH},  with $(\lambda,\mu)=(5,1)$}
\label{tab:TH_new}
\end{table}
\section{Illustrations}
In this section, we present the numerical test in two different situations corresponding to different objectives. We will first consider an academic example with a view to illustrate  the test  itself. Then  we will show how it can be used  to  design with ECP-spaces.

%
\subsection{Academic example}
\label{THTH}
Here, we work with four five-dimensional section-spaces, that is, with $n=4$ and $q=3$. The two section-spaces $\E_0$ and $\E_2$ are the trigonometric spaces spanned on $[t_0, t_1]$ and $[t_2, t_3]$ by the five functions $1,x,x^2,\cos x, \sin x$, while  $\E_1$ and $\E_3$ are the hyperbolic spaces spanned on $[t_1, t_2]$ and $[t_3, t_4]$ by the  functions $1,x,x^2,\cosh x, \sinh x$.  Moreover, in $\E$, the three connection matrices $M_1, M_2, M_3$ are the identity matrices of order five. Accordingly, the space $\E$ is a W-space on $[a,b]=[t_0, t_4]$. For the space $\E$ to be a PEC-space on $\on$
we  have to require that, for $k=0$ and $k=2$, the length $t_{k+1}-t_k$ be less than the so-called {\it critical length} of  five-dimensional trigonometric spaces, which is approximately 8.98, see \cite{CMP}. In that case, for short, we refer to $\E$ as a THTH PEC-space on $\on$. In order to determine if $\E$ is an EC-space on $[a,b]$, we actually apply the numerical test  to $\E$, with
$$
t_1-t_0=\mu, \quad t_2-t_1=t_4-t_3=\lambda, \quad t_3-t_2=5.
$$
Fig.~\ref{fig:gfd_TH} shows the results of the test depending on  the parameters $\lambda$, $\mu$, sampled at equally spaced values in the intervals $\lambda\in [0.05,6.5]$ and $\mu\in[0.05,3.5]$ with sampling step $0.05$. Each point $(\lambda,\mu)$ is depicted in green if the algorithm ends successfully: in that case,  experimentally the space $\E$ is an EC-space on $[a,b]$. The white region (including the boundary of the green region) corresponds to points $(\lambda,\mu)$ where the algorithm stops in either step of T0.  
 In the pink region, the numerical test stops after one iteration (that is, at T1, with $p=1$).

To more precisely illustrate our numerical test, Table \ref{tab:TH_new} concerns the point $(\lambda,\mu)=(5,1)$, located in the pink region. It lists the coefficients of the local expansions of global Bernstein-like bases: on the left, in $\E$, where  the basis $(V_0,\ldots, V_4)$ is given by (\ref{Hermite_endcond1}) and (\ref{Hermite_endcond2});
on the right, in the four-dimensional PEC-space $DL_0\E$ built by generalised differentiation.  In the left part  of Table \ref{tab:TH_new}, all required coefficients are positive, which means that the test can continue, and  that $DL_0\E$ is a PEC-space on $\on$.  By contrast, in the right part of Table \ref{tab:TH_new}, we can see some negative coefficients. Accordingly, the test stops at this stage: we conclude that, experimentally speaking, $\E$ is not an ECP-space on $\on$, that is, it is not an EC-space on $[a,b]$.
Note that the local Bernstein-like bases which have been used to compute these coefficients all satisfy conditions (\ref{Hermite_endcond1}) and (\ref{Hermite_endcond2}) where $a,b$ are replaced by $t_k, t_{k+1}$, respectively. This choice explains the five coefficients equal to 1 in the left part of Table \ref{tab:TH_new} (for $k=0$ and $k=3$).

\begin{remark}
From  Subsections \ref{test1} and \ref{iter} we know that, if an $(n+1)$-dimensional PEC-space $\E$ on $\on$ is not an ECP-space on $\on$, then the test necessarily stops either at T0 or at T1 for some positive $p\leq n-1$. In practice, in most of the many examples we have tested, with various section-spaces, various kinds of connection matrices, and various dimensions, we could observe that the test either stops at T0, or it successfully continues until the end. It was all the more essential to exhibit an example contradicting this observation. The major interest of this example (which may seem somewhat ill conceived at first sight) is thus to clearly point out that the existence of a global Bernstein-like basis in the PEC-space $\E$ (resp. of a global Bernstein basis in the PEC-space $L_0\E$) is not sufficient for $\E$ to be an ECP-space on $\on$ (resp. for $L_0\E$ to be an ECP-space good for design on $\on$) even under the requirement that all convenient coefficients of its local expansions be positive. It secondarily illustrates the possibility of constructing a global EC-space from EC-section spaces. This example therefore also shows that, in a W-space on $[a,b]$ (containing constants), the presence of a positive Bernstein-like basis (a Bernstein basis) relative to $(a,b)$ is not sufficient  for this space to be an EC-space (good for design) on $[a,b]$. See Theorem \ref{th:ECP_NBLB} and its non-piecewise version.
\end{remark}

\noindent
\subsection{Design with ECP-spaces}
\label{design}
Here, we are interested in designing with ECP-spaces. We thus start with an $(n+1)$-dimensional PEC-space $\E$ on $\on$, and we assume that it contains constants. Such a PEC-space $\E$ depends on a number of shape parameters: the entries of the connection matrices plus,  possibly, parameters provided by  the section-spaces themselves.  We would like to determine whether or not $\E$ is an ECP-space good for design on $\on$. We know that this amounts to determining whether or not the $n$-dimensional PEC-space $D\E$ on $\on$ (in which the connection matrices are obtained by deletion of  the first rows and columns in those of $\E$) is an ECP-space on $\on$. Accordingly, applied to $D\E$, the numerical test can be used to experimentally identify  regions of parameters where $\E$ is an ECP-space good for design on $\on$.  Within such experimental regions, we can then try to study how the shape parameters act on the curves for fixed B\'ezier points. 

\smallskip
We first illustrate our purpose with $n=3$, $q=2$, and equispaced knots $t_{k+1}-t_k=1$, $k=0,1,2$.  Moreover, for $q=0, 1, 2$, the section-space $\E_k$ is the cubic polynomial space on $[t_k, t_{k+1}]$, and we assume that the connection matrices $M_1, M_2$ are of the form: 
\begin{equation}
M_1 =
\begin{pmatrix}
1 & 0 & 0 & 0\\
0 & 1 & 0 & 0\\
0 & \beta & 1 & 0\\
0 & \delta & \e & 1\\
\end{pmatrix},
\qquad
M_2 =
\begin{pmatrix}
1 & 0 & 0 & 0\\
0 & 1 & 0 & 0\\[0.5ex]
0 & \beta & 1 & 0\\[0.5ex]
0 & (\beta\e-\delta) & \e & 1\\
\end{pmatrix},
\label{eq:M1M2}
\end{equation}
where $\beta, \delta, \e$ are free real parameters. Our four-dimensional  PEC-space $\E$ is thus composed of piecewise cubic functions on $\on$, which are $C^1$ on $[a,b]$ and geometrically continuous in the sense of continuity of the Frenet frames of order three and of the first two curvatures.  Our choice (\ref{eq:M1M2}) also guarantees that the PEC-space $\E$ is closed under symmetry.

\smallskip

Subsequently we test four different situations where we only have to handle two free parameters.

\begin{figure}[h]
\renewcommand{\thesubfigure}{}
\centering
\subfigure[$\beta=0$]{\includegraphics[width=4.4cm,trim = 0cm 0cm 0cm 0cm, clip]{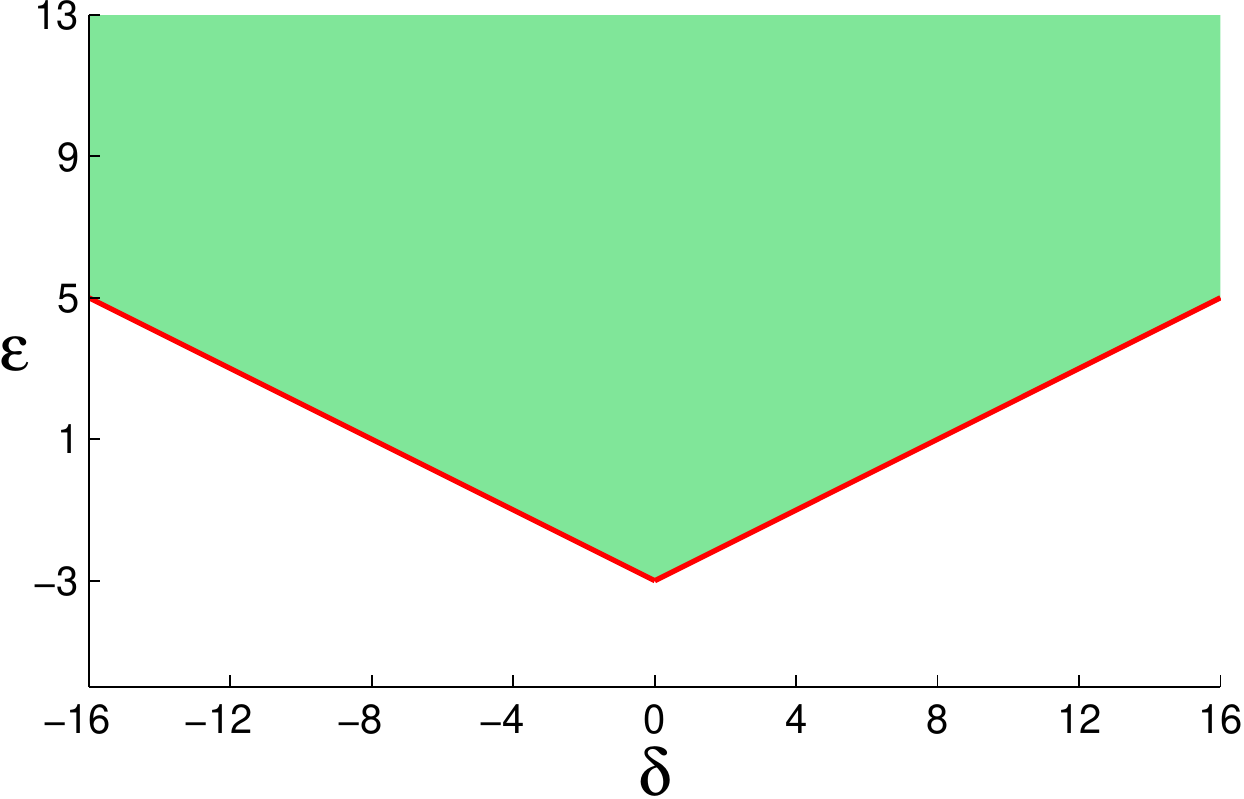}\label{fig:gfd_1}}
\hspace{2cm}
\subfigure[$\varepsilon=0$]{\includegraphics[width=4.4cm,trim = 0cm 0cm 0cm 0cm, clip]{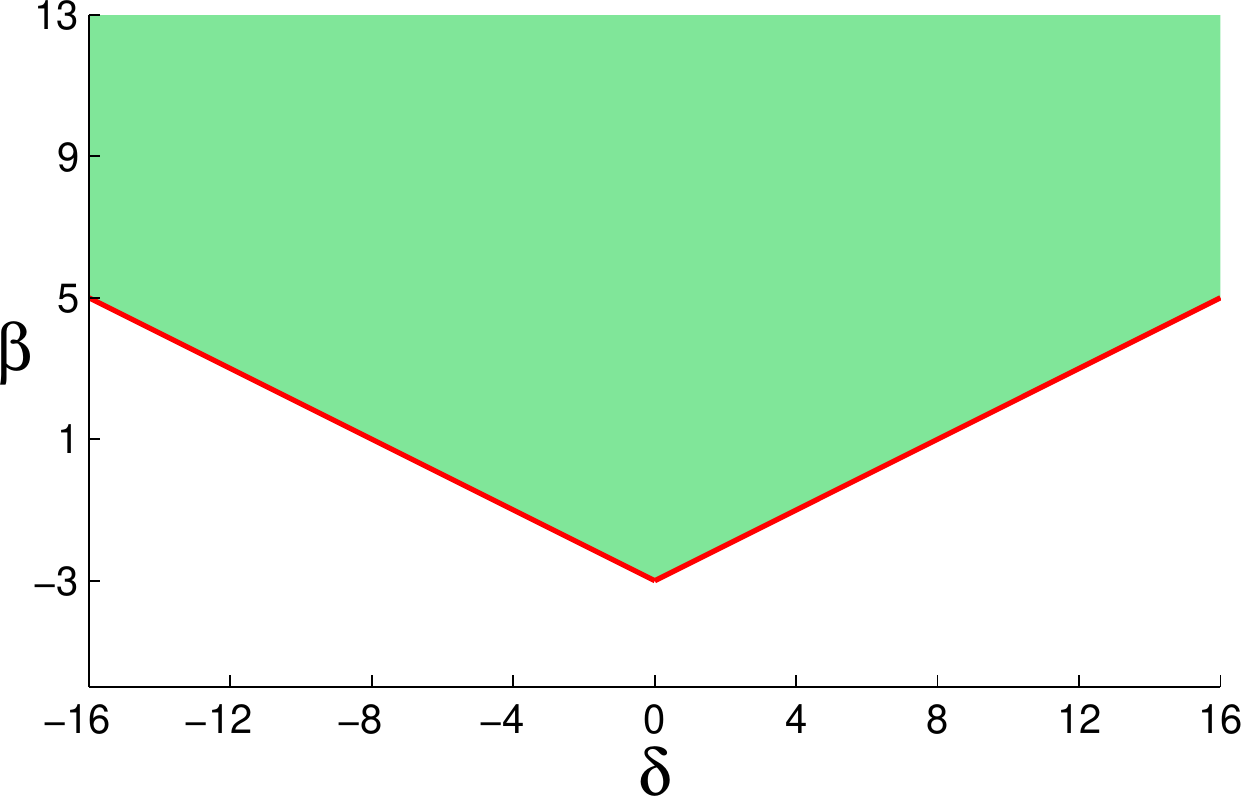}\label{fig:gfd_2}}\\
\subfigure[$\delta=0$]{\includegraphics[width=4.4cm,trim = 0cm 0cm 0cm 0cm, clip]{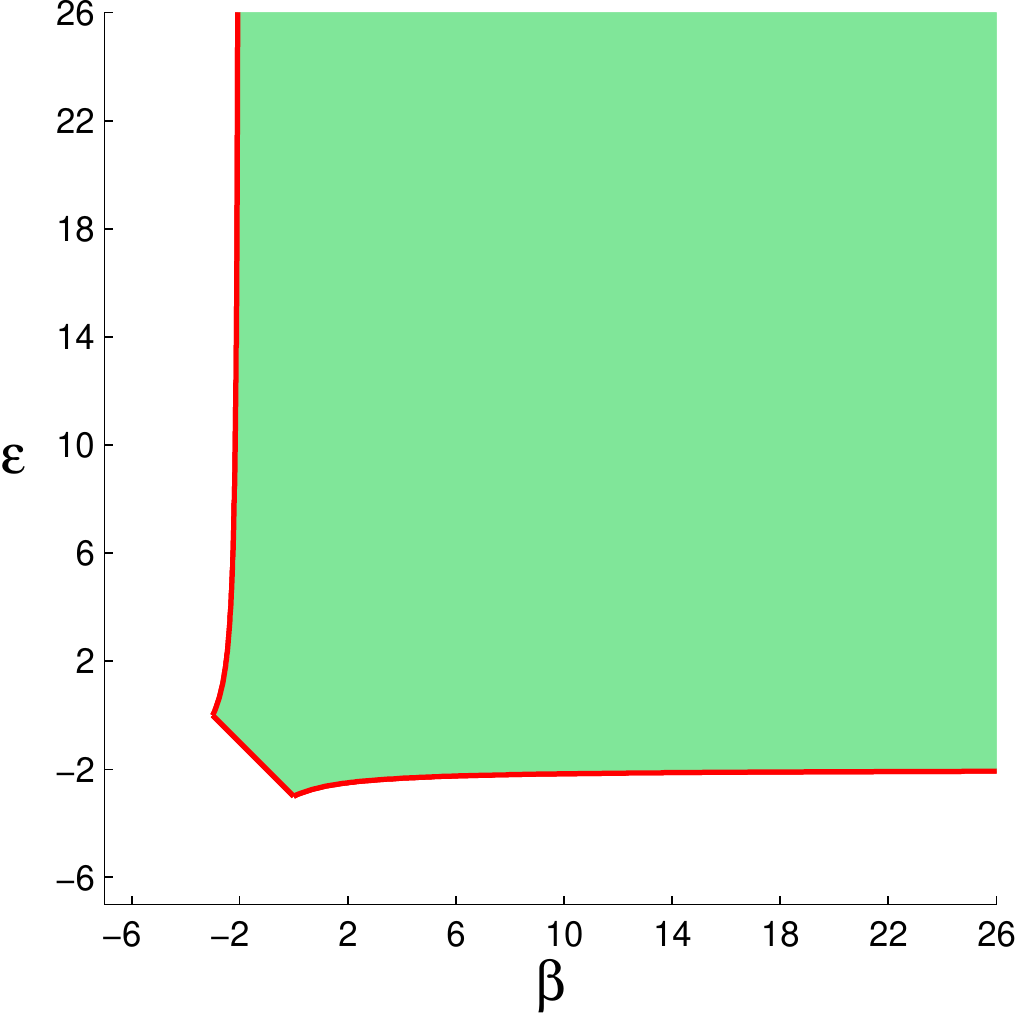}\label{fig:gfd_4}}
\hspace{2.2cm}\subfigure[$\delta=\frac{\beta\varepsilon}{2}$]{\includegraphics[width=4.4cm,trim = 0cm 0cm 0cm 0cm, clip]{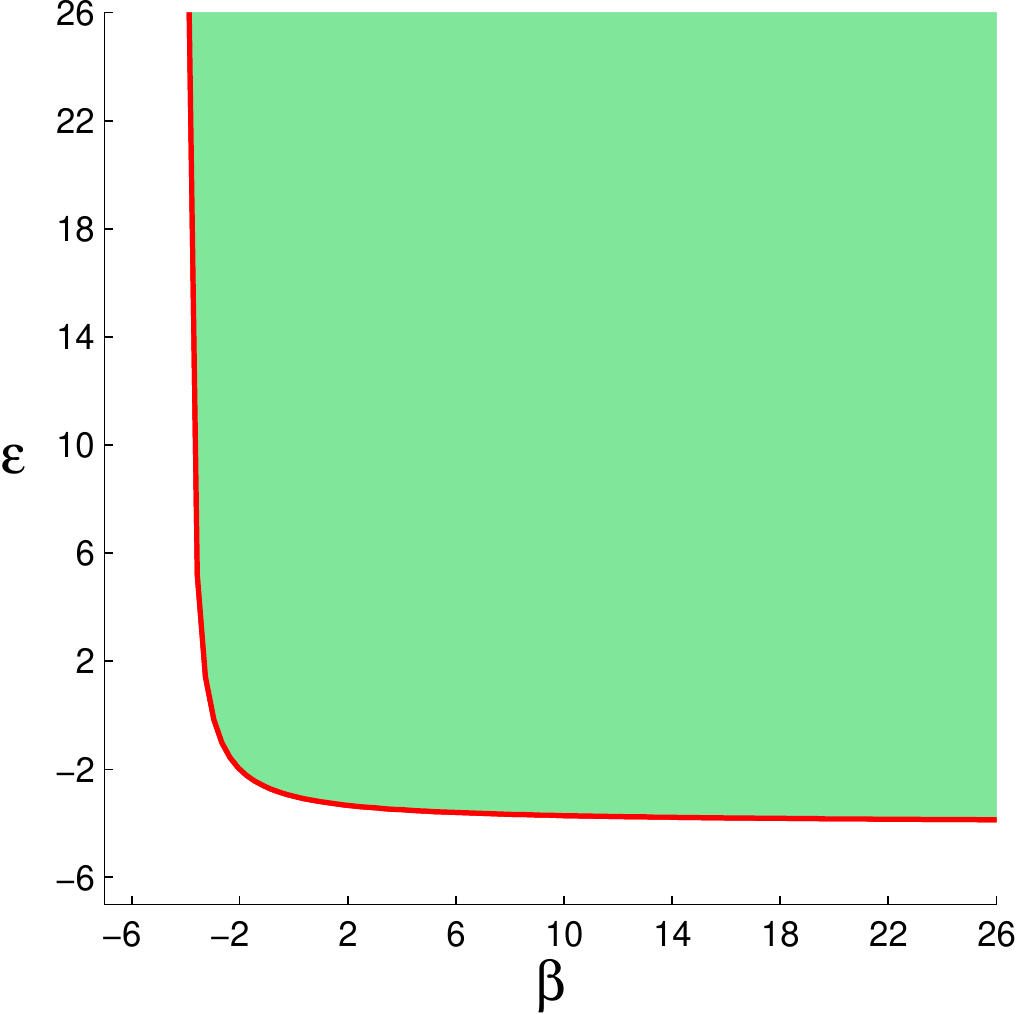}\label{fig:gfd_3}}
\caption{``Good for design" regions (in green) for the PEC-spaces in Subsection \ref{design}.}
\label{fig:gfd}
\end{figure}

\begin{figure}[h]
\centering
\subfigure{\includegraphics[width=0.28\textwidth,trim = 0cm 0cm 0cm 0cm, clip]{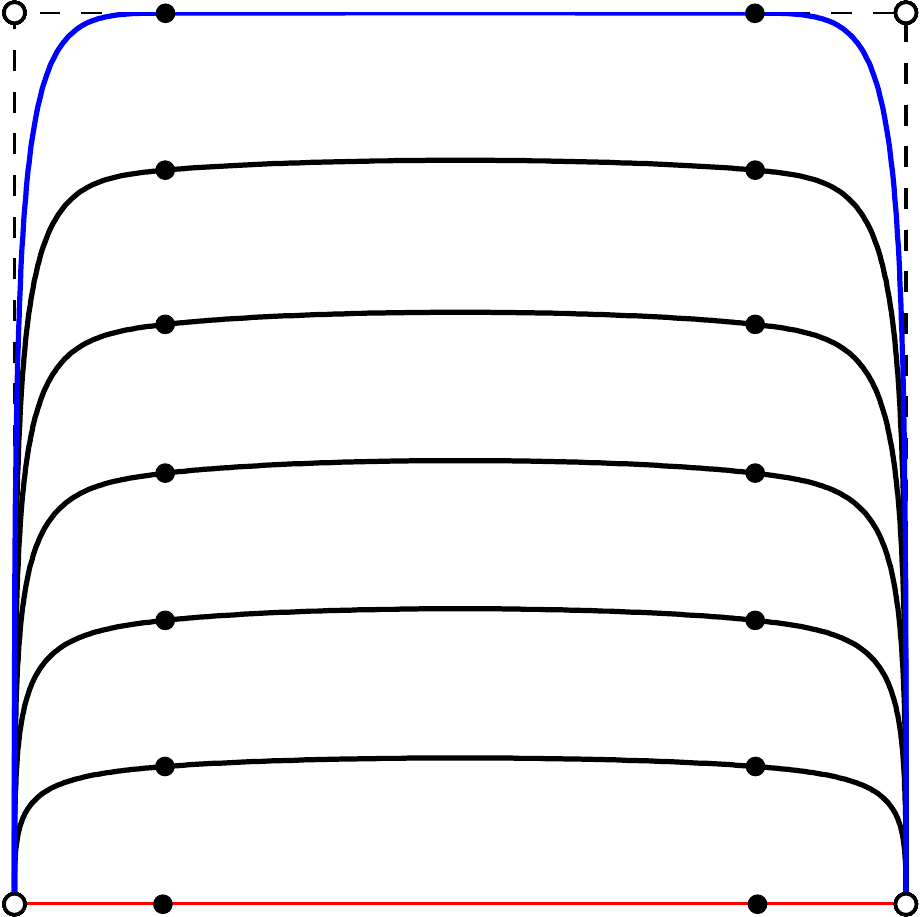}\label{fig_de_b0_1}}\hfill
\subfigure{\includegraphics[width=0.28\textwidth,trim = 0cm 0cm 0cm 0cm, clip]{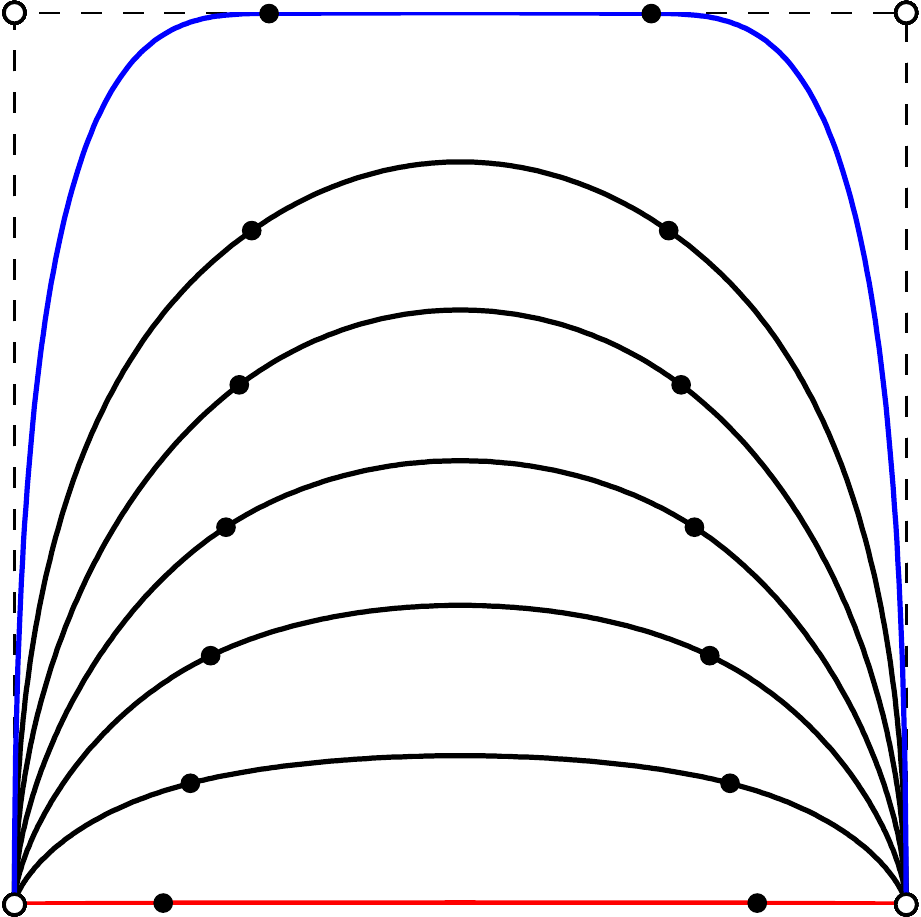}\label{fig_de_b0_2}}\hfill
\subfigure{\includegraphics[width=0.28\textwidth,trim = 0cm 0cm 0cm 0cm, clip]{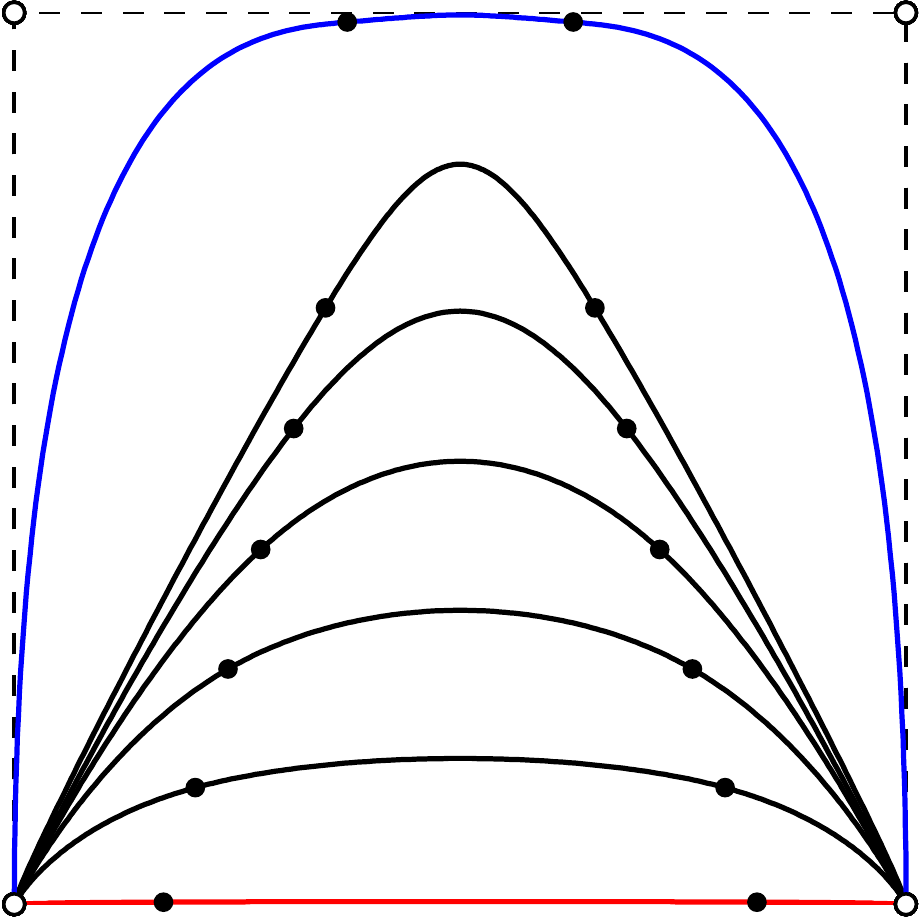}\label{fig_de_b0_3}}\hfill
\caption{Design in the piecewise cubic ECP-space of Subsection \ref{design}, with $\beta = 0$ and $\delta = \alpha(\varepsilon + 3)$.
{\bf Left:} $\alpha = -1.96$, and from bottom to top, $\varepsilon = -2.9974$ (red); $1.6; 4.8; 8.4; 13.5; 24; 1000$ (blue).
{\bf Middle:} $\alpha=0$, $\varepsilon=-2.9974$ (red); $-2.7;-2.3;-1.75;-0.8;1.35;1000$ (blue).
{\bf Right:} $\alpha = 1.96$, $\varepsilon=-2.9974$ (red); $-2.84;-2.59;-2.15;-1.17;2.6;4500$ (blue).}
\label{fig:de_b0}
\vspace{0.5cm}
\centering
\subfigure{\includegraphics[width=0.28\textwidth,trim = 0cm 0cm 0cm 0cm, clip]{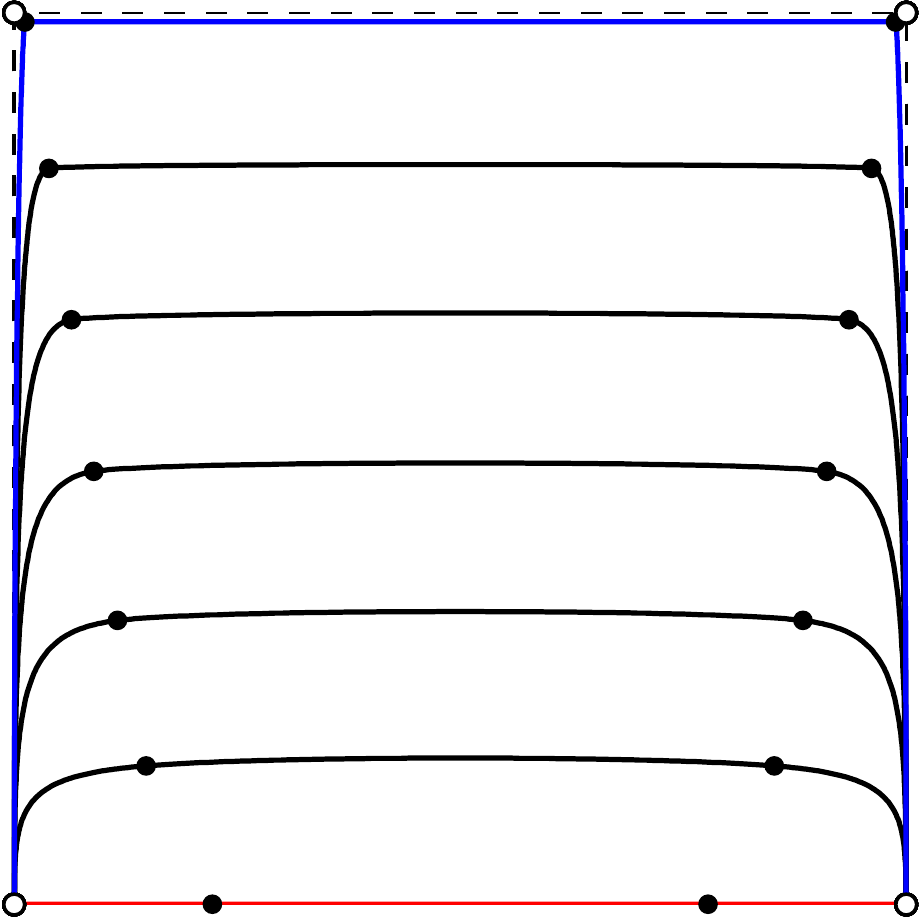}\label{fig_bd_e0_1}}\hfill
\subfigure{\includegraphics[width=0.28\textwidth,trim = 0cm 0cm 0cm 0cm, clip]{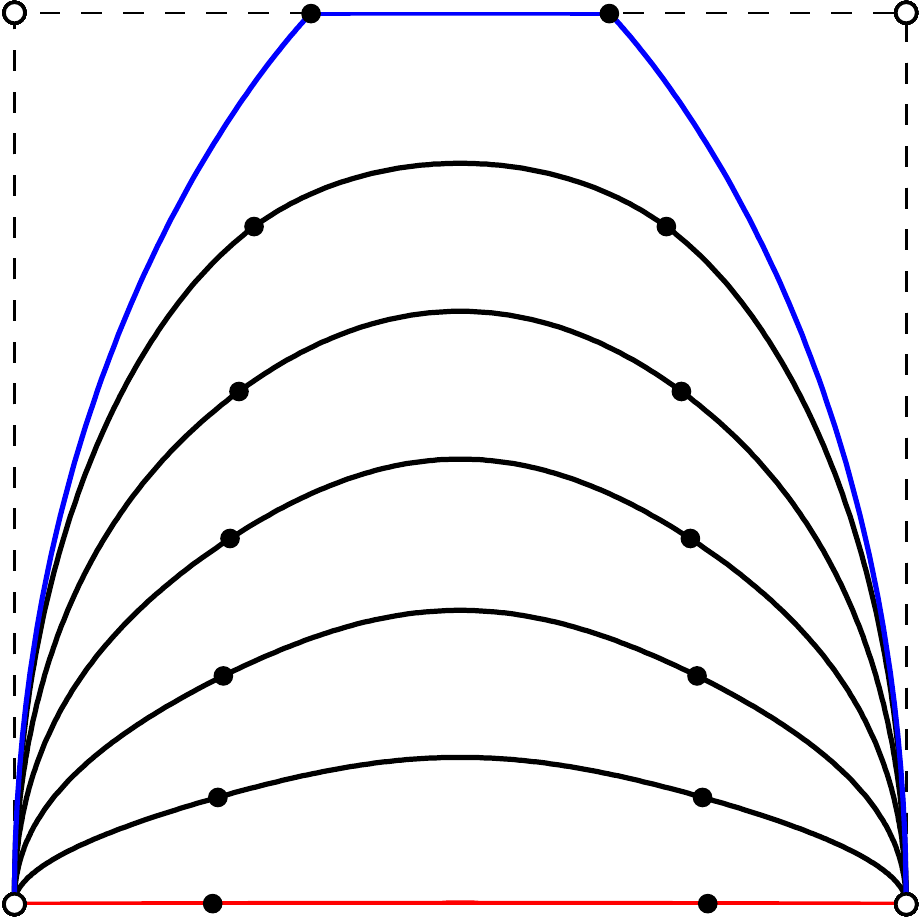}\label{fig_bd_e0_2}}\hfill
\subfigure{\includegraphics[width=0.28\textwidth,trim = 0cm 0cm 0cm 0cm, clip]{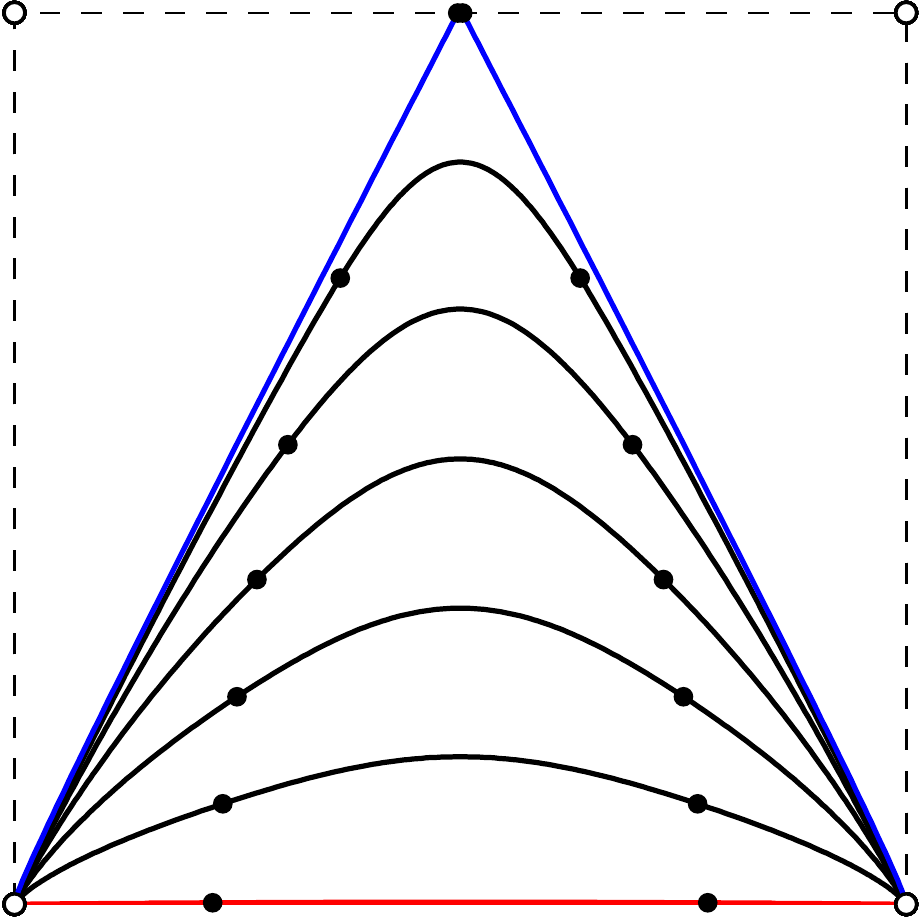}\label{fig_bd_e0_3}}\hfill
\caption{Design in the piecewise cubic ECP-space of Subsection \ref{design}, with $\varepsilon = 0$ and $\delta = \alpha(\beta + 3)$.
{\bf Left:} $\alpha = -1.96$, and from bottom to top, $\beta = -2.9974$ (red); $1.8; 6; 12; 23; 53.5; 1000$ (blue).
{\bf Middle:} $\alpha=0$, $\beta=-2.9974$ (red); $-2.66;-2.26;-1.7;-0.8;1.35;1000$ (blue).
{\bf Right:} $\alpha = 1.96$, $\beta=-2.9974$ (red); $-2.8; -2.5; -2; -1; 1.9; 4500$ (blue).}
\label{fig:bd_e0}
\vspace{0.5cm}
\centering
\subfigure{\includegraphics[width=0.28\textwidth,trim = 0cm 0.0cm 0cm 0.0cm, clip]{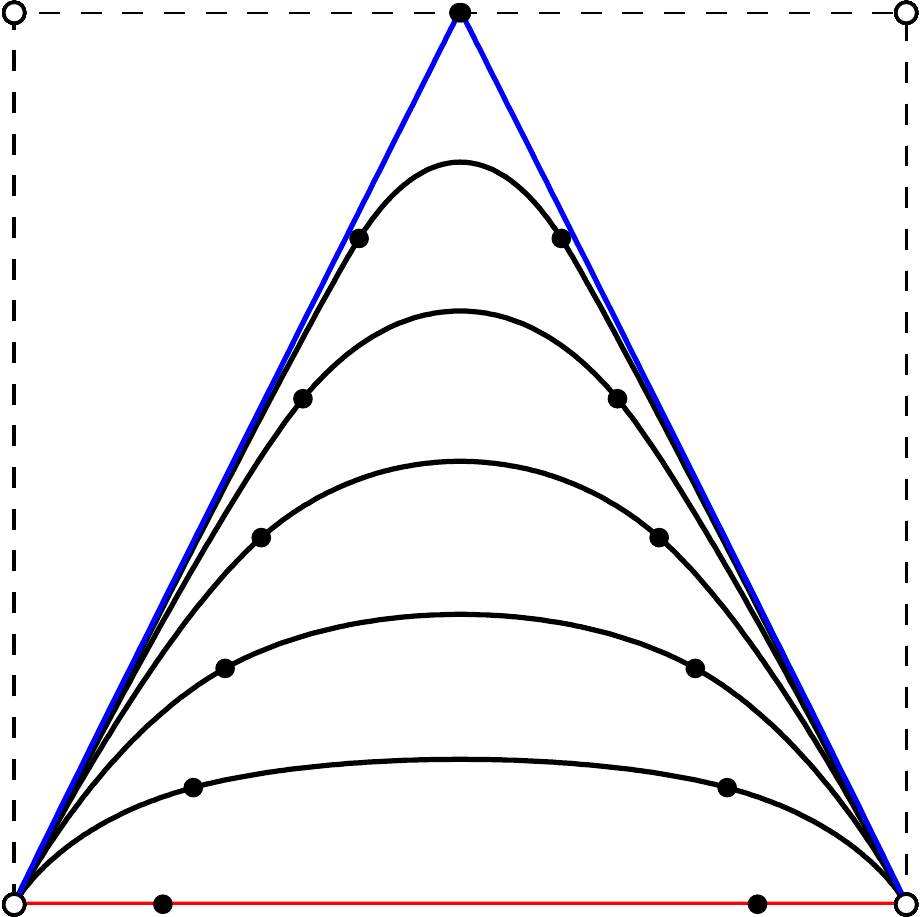}\label{fig:be_d0_1}}\hfill
\subfigure{\includegraphics[width=0.28\textwidth,trim = 0cm 0.0cm 0cm 0.0cm, clip]{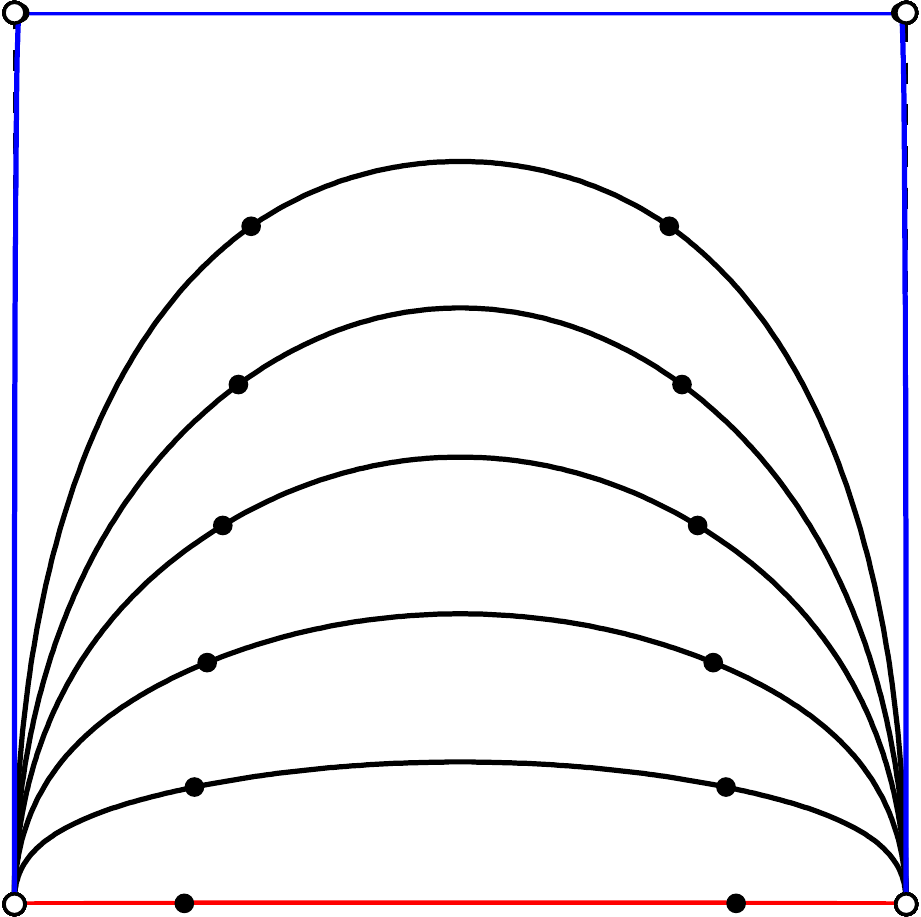}\label{fig:be_d0_2}}\hfill
\subfigure{\includegraphics[width=0.28\textwidth,trim = 0cm 0.0cm 0cm 0.0cm, clip]{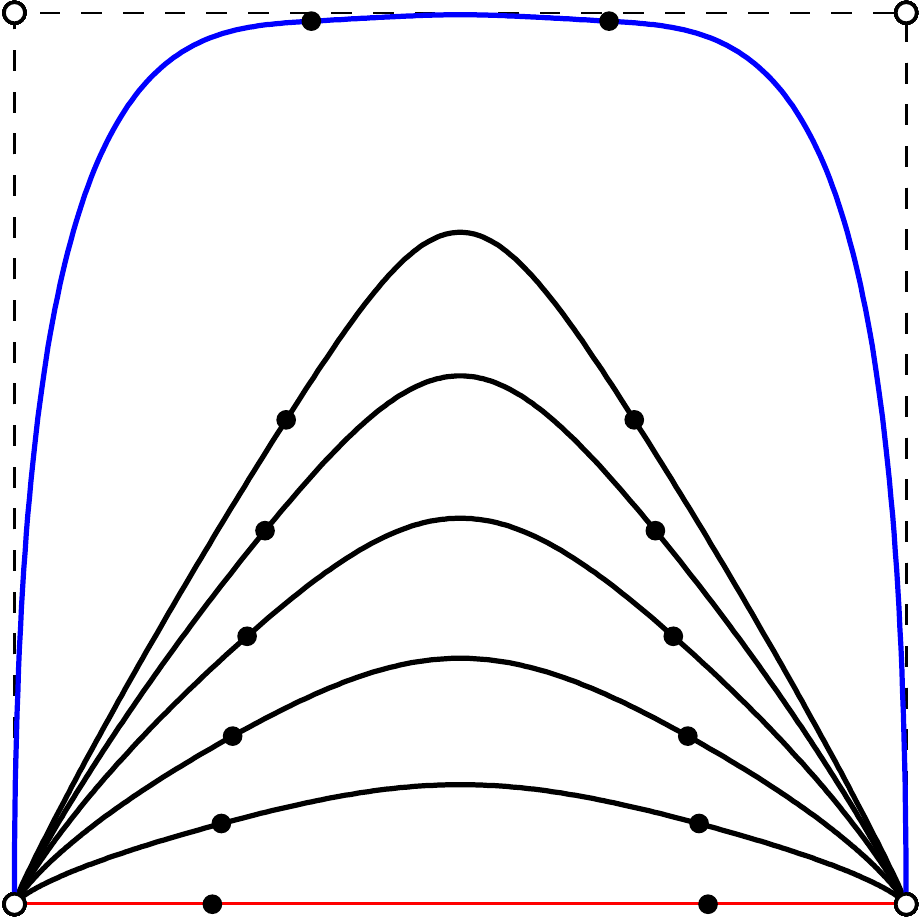}\label{fig:be_d0_3}}\hfill
\caption{Design in the piecewise cubic ECP-space of Subsection \ref{design}, with $\delta = 0$ and
{\bf Left:} $(\beta+1.99)(\varepsilon+1.99)=-2$, with $\beta=-0.003$ (red); $0.21; 0.55; 1.2; 2.6; 7; 20000$ (blue).
{\bf Middle:} $\beta=\varepsilon=-1.499$ (red); $-1.38 ; -1.2 ;-0.9; -0.4; 0.65; 1000$ (blue).
{\bf Right:} $(\beta+1.99)(\varepsilon+1.99)=-2$, with $\varepsilon=-0.003$ (red); $0.21; 0.55; 1.2; 2.6; 7; 20000$ (blue).}
\label{fig:be_d0}
\end{figure}

\smallskip
\noindent
{\bf Case (I):} Here we assume that $\beta=0$. We are thus dealing with $C^2$ piecewise  cubic functions.  The experimental domain of parameters $(\delta, \e)$ where the space $\E$ is an ECP-space good for design on $\on$ is shown in 
Fig.~\ref{fig:gfd}, top, left.  Experimentally, it can be described  by $\delta=\alpha(\e+3)$ with $-2<\alpha<2$ and $\e>-3$. 
The shape effects obtained by varying $(\delta, \e)$ within the experimental green region are illustrated in Fig.~\ref{fig:de_b0}. For a given $-2<\alpha<2$, we can see how the curve changes on the line $\delta=\alpha(\e+3)$ when $\e$ ranges over $]-3, +\infty[$. It moves continuously from the red limit curve obtained for $\e=-3^+$ to the blue one obtained when $\e\rightarrow+\infty$. The values of $\alpha$ which are presented are $\alpha= -1.96$ (left); 0 (middle); 1.96 (right).  The sequence of subfigures efficiently  indicates how the shape effects progressively range from  left to right  as $\alpha$ goes from $-2^+$ to $2^-$.
To better point out  the shape effects,  it is deliberate that we have not used the same values of $\e$ in the three subfigures,
 because the speed of deformation  varies depending on $\alpha$.

\smallskip
\noindent
{\bf Case (II):} Here we assume that $\e=0$. The experimental domain of parameters $(\delta, \beta)$ where the space $\E$ is an ECP-space good for design on $\on$ is shown in 
Fig.~\ref{fig:gfd}, top, right.  Experimentally, it can be described  by $\delta=\alpha(\beta+3)$ with $-2<\alpha<2$ and $\beta>-3$. The shape effects obtained by varying $(\delta, \beta)$ can be observed in Fig.~\ref{fig:bd_e0}, with, similarly to Case (I), from left to right  $\alpha= -1.96; 0 ; 1.96$. Here too the values of the parameter $\beta$ are selected differently depending on $\alpha$ so as  to achieve sufficiently significant variations in shape. 

\smallskip
\noindent
{\bf Case (III):} Here we assume that $\delta=0$.  As shown in Fig.~\ref{fig:gfd}, bottom, left,  the numerical test suggests that the domain of parameters where the space $\E$ is an ECP-space good for design is described  by $(\beta+2)(\e+2)>-2$ and $\beta+\e+3>0$.  In Fig.~\ref{fig:be_d0},  we show the shape effects along the two branches of hyperbolas (the horizontal one in the left picture and the vertical one in the right picture). For numerical reasons, we have to remain sufficiently far from the asymptotes. This is why in both cases the selected pairs $(\beta, \e)$ are taken from the hyperbola
$(\beta+1.99)(\e+1.99)=-2$. The curves in the middle picture represent pairs on the diagonal, that is,  $\e=\beta >-1.5$. As can logically be expected from the previous pictures, along the segment $\beta+\e+3=0$ (for instance, on the segment $\beta+\e+2.9=0$), there  is hardly any change in the curves which all nearly coincide with the segment joining the extreme B\'ezier points. %

\begin{figure}
\centering
\subfigure{\includegraphics[width=0.28\textwidth,trim = 0cm 0.0cm 0cm 0.0cm, clip]{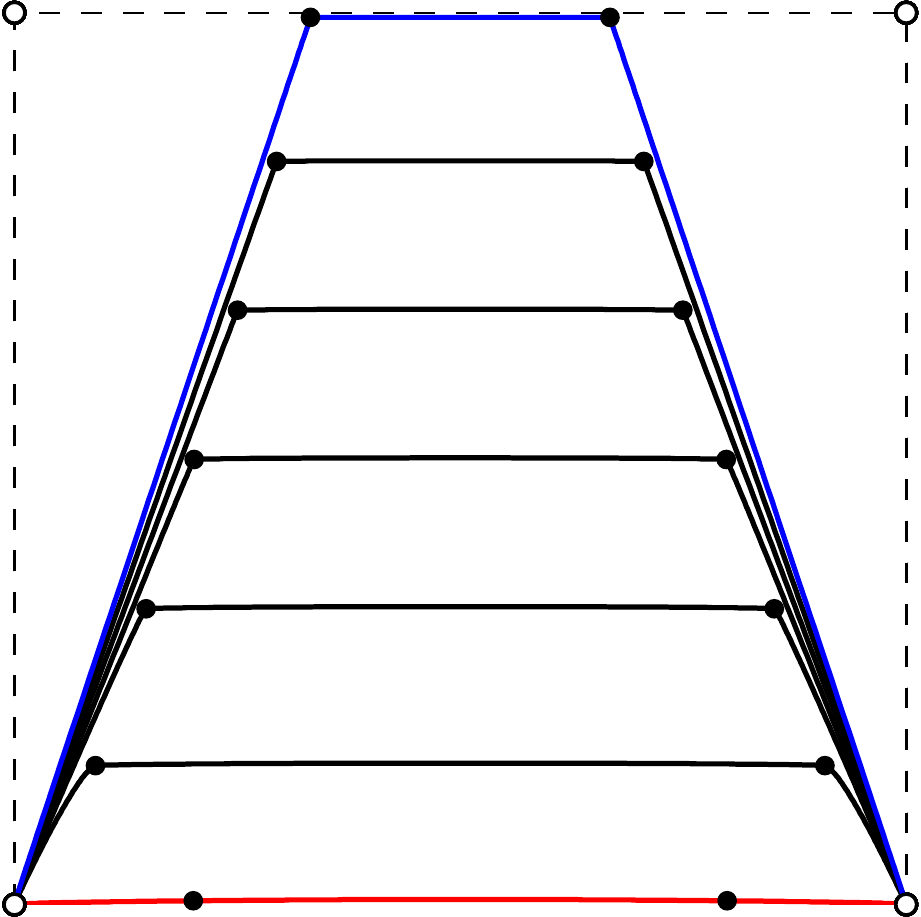}\label{fig:be_d05be_1}}\hfill
\subfigure{\includegraphics[width=0.28\textwidth,trim = 0cm 0.0cm 0cm 0.0cm, clip]{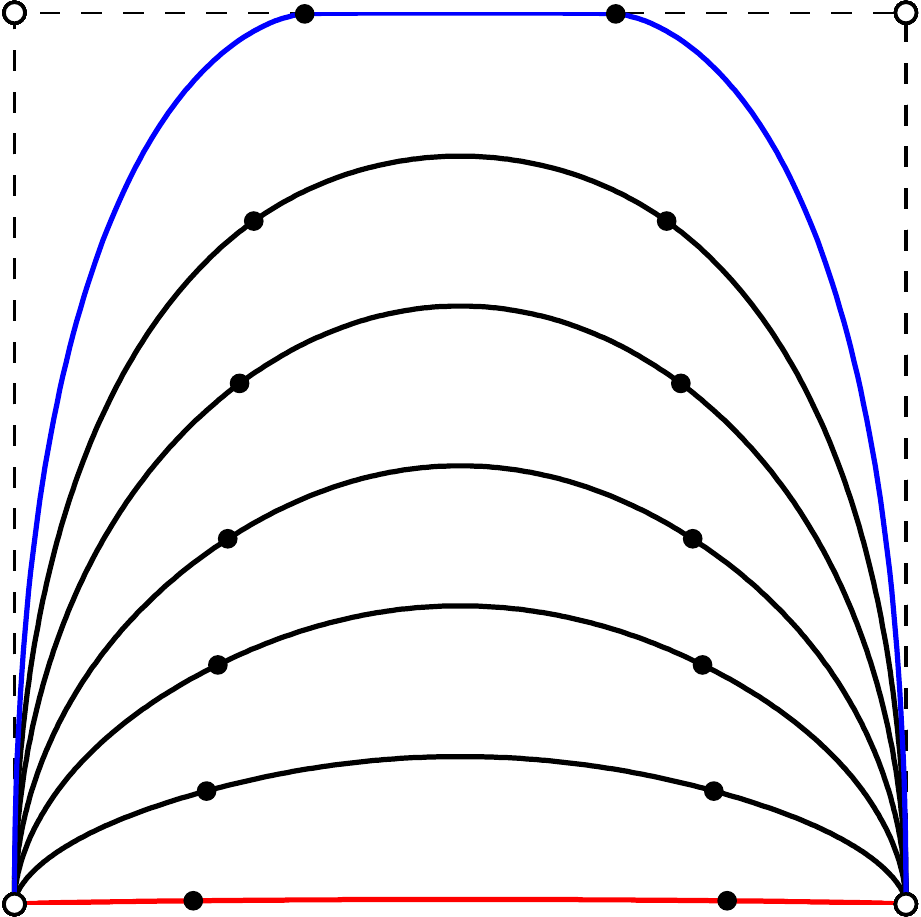}\label{fig:be_d05be_2}}\hfill
\subfigure{\includegraphics[width=0.28\textwidth,trim = 0cm 0.0cm 0cm 0.0cm, clip]{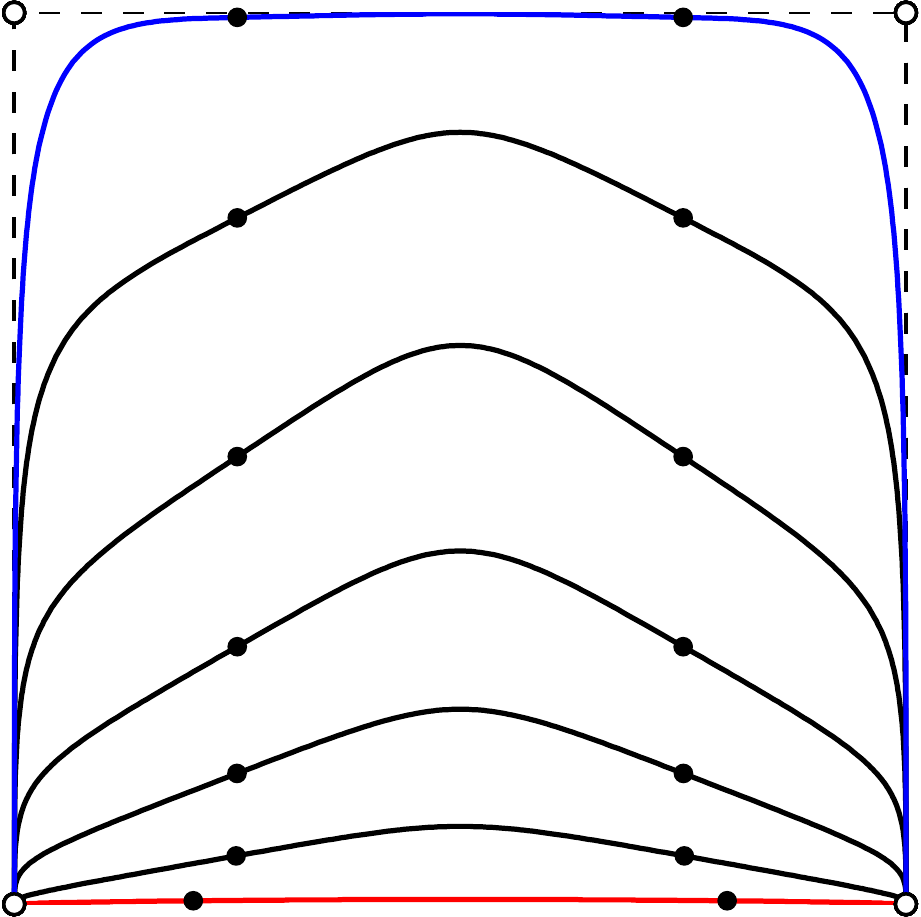}\label{fig:be_d05be_3}}\hfill
\caption{Design in the piecewise cubic ECP-space of Subsection \ref{design}, with  $d=\frac{be}{2}$.
{\bf Left:} $(\beta + 3.99)(\varepsilon + 3.99) = 4$, with
 $\beta = -1.99$ (red); $75; 200; 400; 800; 2000; 75000$ (blue).
{\bf Middle:} $b=e=-1.99$ (red); $-1.7; -1.37; -1;-0.4; 0.7; 50$ (blue).
{\bf Right:} $(\beta + 3.96)(\varepsilon + 3.96) = 4$, with
 $\varepsilon = -1.99$ (red); $75; 200; 400; 800; 2000; 75000$ (blue).
}
\label{fig:be_d05be}
\vspace{0.5cm}
\centering
\subfigure{\includegraphics[width=0.28\textwidth,trim = 0cm 0.0cm 0cm 0.0cm, clip]{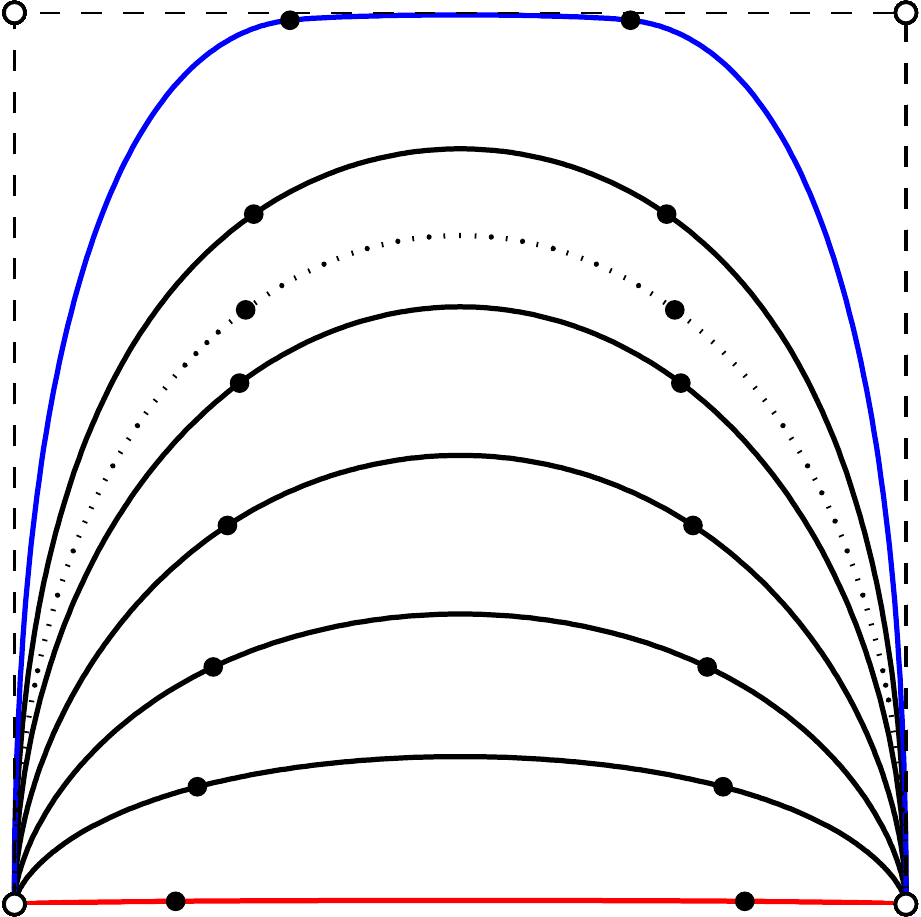}\label{fig:G3_1}}\hfill
\subfigure{\includegraphics[width=0.28\textwidth,trim = 0cm 0.0cm 0cm 0.0cm, clip]{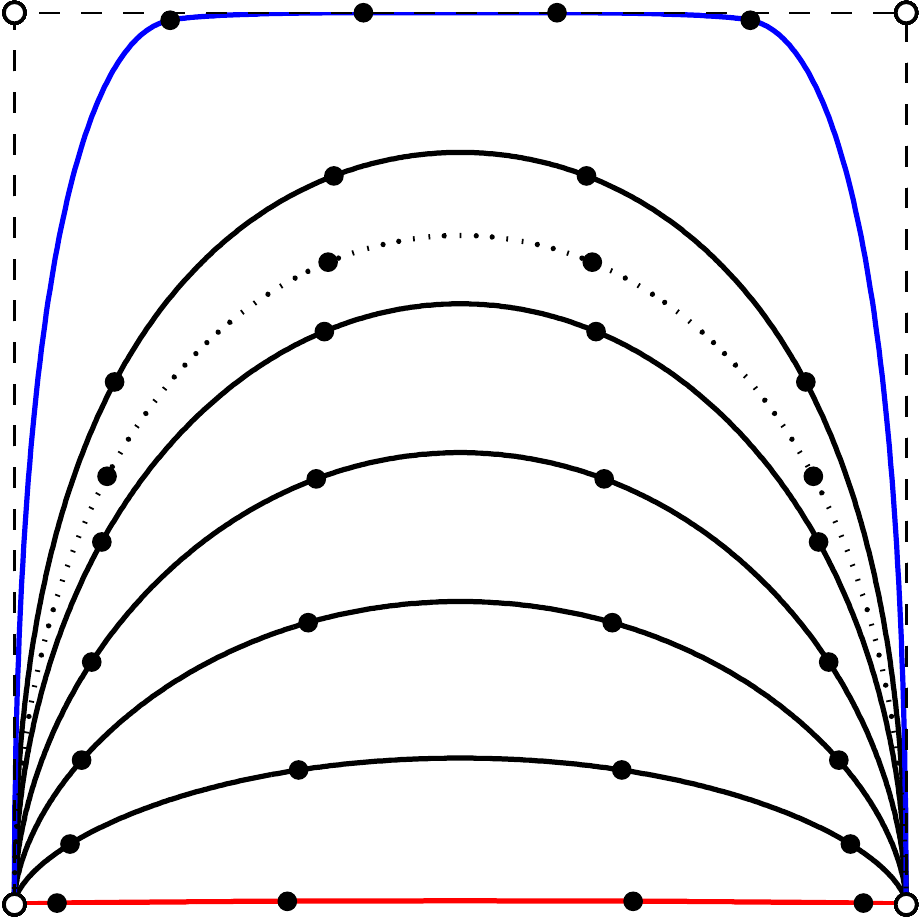}\label{fig:G3_2}}\hfill
\subfigure{\includegraphics[width=0.28\textwidth,trim = 0cm 0.0cm 0cm 0.0cm, clip]{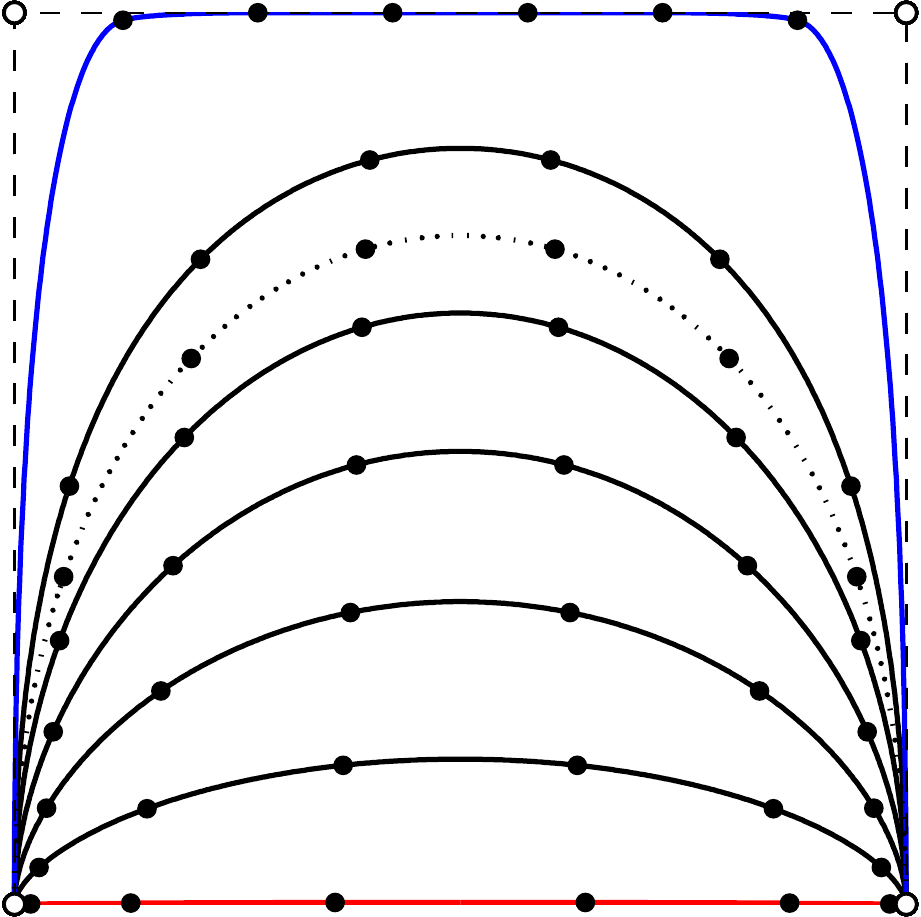}\label{fig:G3_3}}\hfill
\caption{$G^3$ piecewise cubics with 3, 5, 7 sections, depending on the parameter $\beta$, with everywhere
the same matrix.
{\bf Left:} $\beta=-0.9$ (red); $-0.79; -0.66; -0.47; -0.2; 0.4; 10$ (blue).
{\bf Middle:} $\beta=-0.37$ (red); $-0.315; -0.25; -0.175; -0.07; 0.13; 10$ (blue).
{\bf Right:} $\beta = -0.195$; (red) $-0.165;-0.13; -0.09; -0.04; 0.07;10$ (blue).
}
\label{fig:G3}
\end{figure}

\smallskip
\noindent
{\bf Case (IV):} Here we assume that $\delta=\frac{\beta\e}{2}$, or equivalently, that the matrices $M_1$ and $M_2$ are equal. As shown in Fig.~\ref{fig:gfd}, bottom, right,  the numerical test suggests that the domain of parameters where the space $\E$ is an ECP-space good for design is described  by $(\beta+4)(\e+4)>4$ and $\beta+4>0$.  In Fig.~\ref{fig:be_d05be} we show the shape effects along the boundary of the domain.  As in the previous case, for numerical reasons,  the corresponding pairs $(\beta, \e)$ are selected on the hyperbola
$(\beta+3.96)(\e+3.96)=4$, the left picture corresponding to its vertical half, and the right one to the horizontal half. \; The curves

\medskip
\clearpage
\noindent
in the middle picture represent pairs on the diagonal, namely $\e=\beta>-2$.

 In Fig.~\ref{fig:G3}, the four-dimensional ECP-space $\E$ is again composed of piecewise cubics, with, from left to right $q=2$; 4 ; 6 interior knots.  We take everywhere the connection matrix  already used  in Case (IV), but we additionally assume that $\e=3\beta$. We thus have $G^3$ piecewise cubics (that is, $C^3$ with respect to the arc-length) depending on only one parameter $\beta$.  The test tells us where to choose this parameter depending on the number of sections. For comparison, the dotted curve represents  the ordinary cubic, obtained here with $\beta=0$.  We can see both the efficiency of the parameter $\beta$, how important it is for shape effects to allow negative values of the parameters, along with the evolution when the number of section increases.

%

\section{Concluding comments}

We have presented a numerical test to determine whether a given PEC-space on $\on$ is an ECP-space on $\on$, whose main application  is to determine whether or not a PEC-space on $\on$ which contains constants  is an ECP-space good for design on $\on$. The example of piecewise cubic PEC-spaces considered in Subsection \ref{design} will serve as a basis for commenting on the usefulness and the limits of this test.


\subsection{Shape effects from ECP-spaces}
\medskip
The usefulness of the test is mainly justified by the practical interest of ECP-spaces for design. Indeed, compared to polynomial spaces, and even to EC-spaces, the interest of using ECP-spaces lies in the fact that they mix within the same space two different families of  parameters: those provided by their section-spaces (whose effects are known in the most classical examples)
along with  those provided by the connection matrices (some  effects  of which can be foreseen, \eg changes in the curvatures).  The interactions between these two families can be used to create new interesting shape effects.

The expression ``shape parameter" refers to a parameter which can be used to slightly modify/improve the shape of a curve while  keeping its general aspect, determined by fixed control points.  This is expected to be an interactive process, and the designer should therefore be able to predict how changing the parameter will affect the shape of the curve. In this regard, it is not reasonable to deal with a great number of parameters at the same time. Now, concerning four-dimensional PEC-spaces on $\on$, each  connection matrix is defined by six parameters, while each section space provides us with  up to 3 free shape parameters. The total amount of parameters can thus reach $6q+3(q+1)=9q+3$. Even with a small number of sections, we cannot expect to efficiently control the shape deformations with such a number of parameters varying simultaneously. This is the reason why, in order to present relevant examples of shape variations, we  first have drastically reduced the number of free connection parameters to three, via various arguments. Additional possibilities could be obtained in the general case: for instance, requiring unit diagonals does not allow changes in the first two curvatures, which could be interesting to produce special shape effects. We also have deliberately  limited ourselves to the simplest piecewise cubic case where no shape parameter comes from the section-spaces themselves, at least after fixing the interior knots.

Even under such limitations, the samples of geometrically continuous piecewise cubic curves presented in Figs.~3 to 7 give clear evidence of the amplitude of the shape effects by comparison with known four-dimensional EC-spaces. For an efficient analysis of these effects, each of these figures  actually involves at most two free parameters. For example, we can go continuously from the segment joining the first and last B\'ezier points to the control polygon itself, and this can even be done in different ways depending on the parameters used (see Fig.~3 and 4, Left,  compared to Fig.~5, Middle).
In each case, the intuitive prediction of the shape effects is made possible by investigating the limit cases, which makes it necessary to first determine numerically the boundary of the ``good for design" region. In our examples we can also observe some redundancy in the shape effects, but our purpose is not to investigate this in more details.

\subsection{Theoretical ECP-spaces versus Numerical}

Quite obviously, whenever possible, it should always be preferred to determine  the exact region in the space of all parameters within which the given PEC-space is an ECP-space for design. Not only is this more satisfying mathematically speaking, but it also helps to understand what are the true shape parameters, and it permits a more efficient use of the limit cases.

\medskip
Assuming that $q=1$, consider any four-dimensional PEC-space on $\on$ containing constants, with given section-spaces, the connection matrix at the unique interior knot $t_1$ being the matrix $M_1$ in (\ref{eq:M1M2}), but now with any diagonal positive entries $\alpha, \gamma, \zeta$. Then, the  exact ``good for design" region in the space $\RR^6$ can always be deduced from Proposition 7.1  of \cite{CAGD2016}.  Note that it also involves  the possible parameters provided by the section-spaces, that is, up to 6 additional parameters. This is an example where using the numerical test would not be so relevant.

By way of illustration, consider piecewise  cubics with $t_2-t_1=t_1-t_0=1$. Then, from Proposition 7.1  of \cite{CAGD2016} we can derive that the ``good for design" region  is defined by the following inequalities
\begin{equation}
\label{two sections}
\e+2(\gamma+\zeta)>0, \quad -2(\beta+\e+\alpha+2\gamma+\zeta)<\delta<\frac{(\e+2\zeta)(\beta+2\alpha)}{\gamma}+2(\alpha+\zeta).
\end{equation}
In the symmetric case (that is, when  $\alpha=\gamma=\zeta=1$ and $\delta=\beta\e/2$), the conditions (\ref{two sections}) reduce to $\beta+4>0$ and $\e+4>0$. Illustrations can be found in \cite{CAGD2016}.

When increasing the number of four-dimensional sections, in principle an explicit exact  description of the ``good for design" region can be achieved by iteration, adding one more section-space at each step, but  the process can be laborious. Since this has not been explicitly done yet, the numerical approach offers an interesting alternative. Obviously, the numerical test  is even more useful for higher dimensional section-spaces, where no explicit conditions have been obtained so far.

From the computational point of view, there is no problem in applying the test with different higher dimensional section-spaces, whatever the number of parameters they involve,
with the most general connection matrices, depending on the interior knots. The limitation is of conceptual nature: as already mentioned, too many parameters do not permit an intuitive handling of the corresponding shape effects. 
\subsection{Piecewise Chebyshevian splines}

As reminded in  the introduction, a  given Piecewise Chebyshevian spline space  is good for design if and only if it is based on an ECP-space good for design. Since the numerical test aims at determining ECP-spaces good for design, it is natural to discuss its possible benefits  in the study of splines.

\medskip
Firstly, we would  like to draw the reader's attention to the fact that
 the numerical test presented in this article is not an appropriate tool  for determining when a  given Piecewise Chebyshevian spline space $\spl$ is good for design.

To illustrate the previous statement, consider   the four-dimensional piecewise cubic PEC-space $\E$ on $\on$, with $q\geq 2$ and knot spacing equal to 1, and   with everywhere the same connection matrix (\ref{eq:M1M2}) obtained with $\delta=\e=0$. Now, for $q=2$, numerically speaking, the space  $\E$ is an ECP-space good for design on $\on$ if and only $\beta>-3$. When increasing the number of sections, the ``good for design" region of the parameter $\beta$ remains of the form $]\beta_0, +\infty[$, but $\beta_0$ increases a lot. This is shown in Fig.~8, Left,  where we can see  the ``good for design" region for $q=2, 4,6, 8$ interior knots. This suggests that $\beta_0$ tends to 0 when $q\rightarrow+\infty$.  As a matter of fact, this is confirmed when applying successively the test with 20, 40, 60 section-spaces: it indicates that $\beta_0$ lies in $]\!-0.0979,-0.0978[$,  $]\!-0.0247,-0.0246[$,
$]\!-0.011,-0.0109[$, respectively.

Let us now denote by $\spl$ the spline space based on $\E$ obtained when each interior knot of $\T$ is simple. In other words, $\spl$ is composed of all $C^1$ piecewise cubic splines $S:[a,b]=[t_0, t_{q+1}]\rightarrow \RR$ satisfying the connection conditions
\begin{equation}
\label{beta}
S''(t_k^+)=\beta S'(t_k) +S''(t_k^-), \quad k=1, \ldots, q.
\end{equation}
From $\E \subset \spl$ we can only deduce that, for $\beta>\beta_0$,  the spline space $\spl$ is good for design.
For comparison, let us now investigate $\spl$ on the theoretical side. We know that $\spl$ is good for design if and only we can find generalised derivatives associated with its section-spaces with respect to which the connection conditions (\ref{beta}) are expressed by the identity matrix \cite{howto}.  While this beautiful theoretical characterisation is difficult to exploit in practice for PEC-spaces, in the spline space $\spl$ it is straightforward to transform it into the very simple inequality $\beta+4>0$, independently of the number of section-spaces, see \cite{beyond, NUMA16}. Accordingly, the theory developed in \cite{howto} tells us that, independently of $q$, given any $\beta>-4$, we can build  infinitely many piecewise cubic ECP-spaces $\overline\E$ good for design on $\on$ so that  $\overline\E\subset \spl$. For $-4<\beta\leq \beta_0$, none of them is equal to $\E$. Without going into details, we can even require  $\overline\E$ to have the same matrix as in Case (IV) everywhere. Combining the theoretical condition $\beta+4>0$ and the numerical results presented in Fig.~8, Right, suggests that the  corresponding ``good for design" regions are defined by $(\beta+4)(\e+4)>\mu$ and $\beta+4>0$, the real number $\mu$ increasing from $4$ to $16^-$ as the number of section-spaces increases from 3 to $+\infty$.

\begin{figure}[h]
\renewcommand{\thesubfigure}{}
\centering
\subfigure[$\delta=\varepsilon=0$]{\includegraphics[height=5.0cm,trim = 0cm 0cm 0cm 0cm, clip]{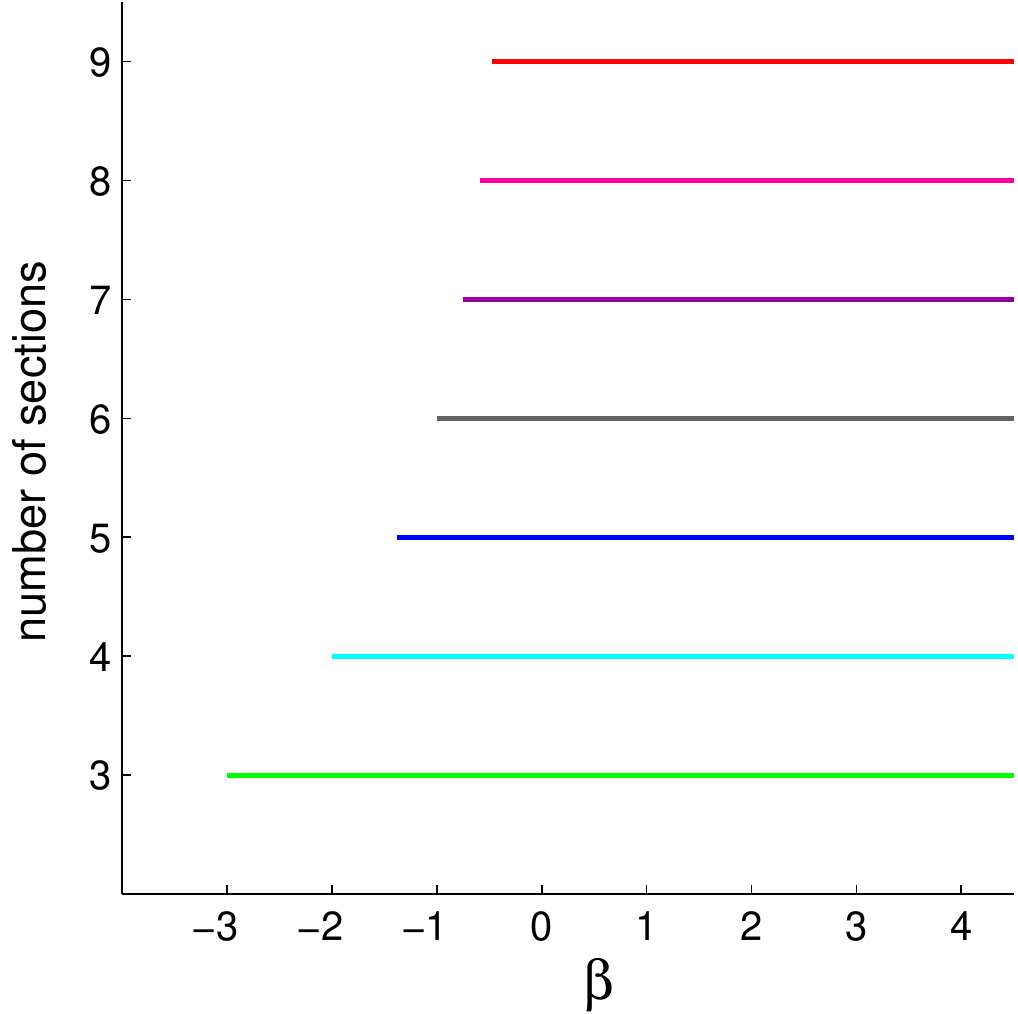}\hfill\label{fig:8_left}}
\hspace{2cm}
\subfigure[$\delta=\frac{\beta\varepsilon}{2}$]{\includegraphics[height=5.0cm,trim = 0cm 0cm 0cm 0cm, clip]{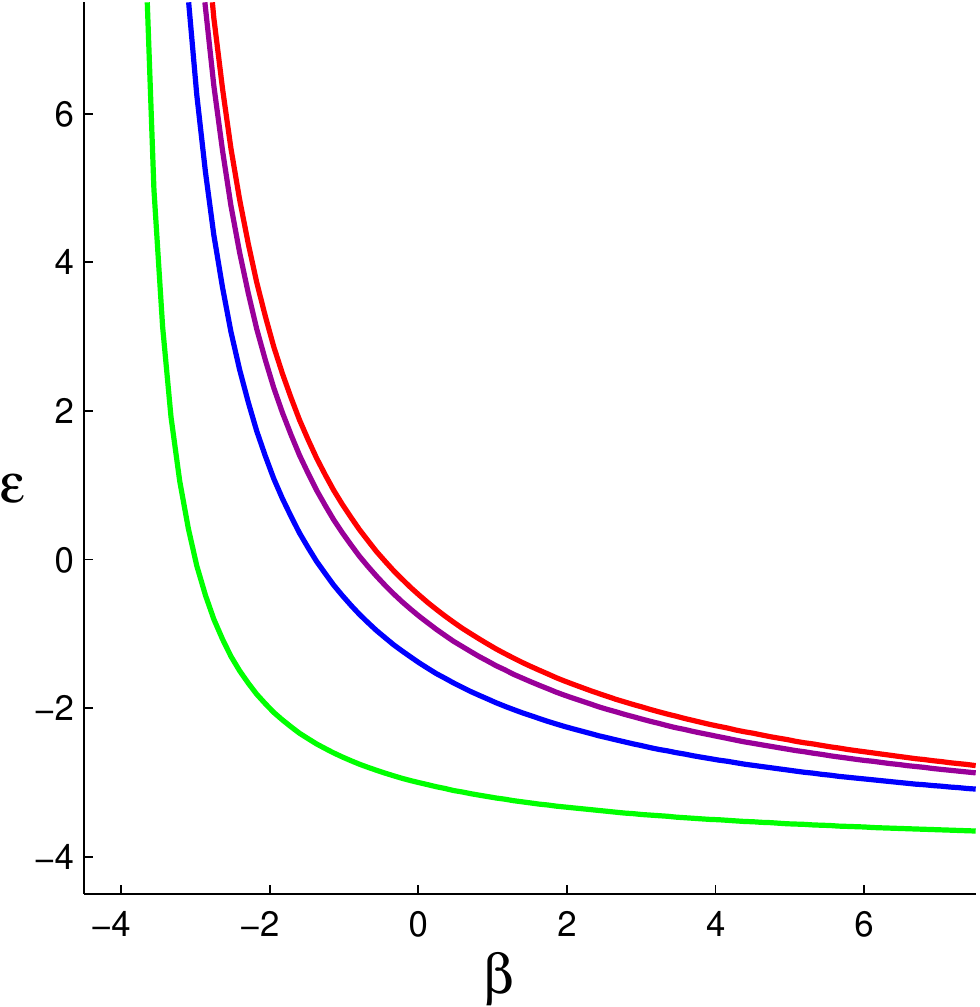}\hfill\label{fig:8_right}}
\caption{``Good for design regions" with increasing number of section-spaces and the same connection matrix everywhere. {\bf Left:} $\delta=\e=0$. {\bf Right: $\delta=\beta\e/2$}, with $3, 5, 7, 9$ section-spaces. }
\label{fig:figure8}
\end{figure}

In spite of the previous limitation, the numerical test can  efficiently  contributes to the spline context. We shall now briefly explain why and how.

Considering splines with non-negative rather than positive interior multiplicities is not usual in CAGD. Still, in the context of Piecewise Chebyshevian splines, it is absolutely natural, and this already proved to efficiently locally increase the flexibility of the spline curves and therefore the shape possibilities \cite{PJMLV}.

Clearly, the class of 	all Piecewise Chebyshevian spline spaces good for design with non-negative multiplicities coincides with the class of all spline spaces good for design with only positive multiplicities and with ECP-spaces good for design as section-spaces.
Once a given PEC-space on $\on$ containing constants is known to be an ECP-space good for design on $\on$, it can therefore serve as one section-space on the interval $[t_0, t_{q+1}]$, with a view to building Piecewise Chebyshevian spline spaces good for design  with positive multiplicities at the  knots $t_0$ and $t_{q+1}$, now possibly considered as interior knots. We can then adapt to them necessary and sufficient conditions for spline spaces with positive multiplicities to be good for design, if such conditions are available, which is the case  for four- and five-dimensional section-spaces (see, for instance, \cite{beyond, NUMA16}). We will show how to more generally exploit this in a further paper, the numerical test being particularly suitable for such applications.

\subsection{Extensions}

The numerical
test presented here can easily be adapted to the larger context of PQEC-spaces (Piecewise Quasi Extended Chebyshev spaces)  on $\on$ obtained when the section-spaces are only assumed to be QEC-spaces (Quasi Extended Chebyshev spaces) on their intervals. The only two differences are as follows: firstly, the connections  must be expressed in a slightly different way; secondly, Bernstein(-like) bases have to be replaced by Quasi Bernstein(-like) bases (see \cite{CAGD2016} and references therein). These differences are due to the fact that the section-spaces do not guarantee unisolvent  Taylor interpolation problems.  As in the present paper,  the underlying tools are blossoms and their pseudoaffinity (\cite{CAGD2016} and references therein concerning QEC-spaces).  A spline version of the test can also be developed. Nevertheless, whenever possible, splines should preferably be treated as explained in the second part of Subsection 5.3.

\end{document}